\documentclass{article}
\usepackage{authblk} 
 
\usepackage[utf8]{inputenc} 
\usepackage[USenglish]{babel} 
\usepackage{graphicx} 
\usepackage{comment} 
\usepackage{amsmath}
\usepackage{amsfonts}
\usepackage{amssymb}\usepackage{amsmath} 
\usepackage{amssymb}
\usepackage[singlespacing]{setspace}
\usepackage{mathtools} 
\usepackage{stmaryrd}
\usepackage{romanbar}
\usepackage[table]{xcolor}
\usepackage{booktabs}
\usepackage{multirow}

\usepackage{amsthm}
\usepackage{abstract} 
\usepackage{float}
\usepackage[mathscr]{euscript}
\usepackage{caption}
\usepackage[singlelinecheck=false,justification=RaggedRight,format=plain]{caption}
\usepackage{subcaption}
\usepackage[inner=2.50cm, outer=2.50cm, top=2.50cm, bottom=3.50cm]{geometry}
\usepackage{mathabx}
\usepackage{tikz-cd}
\usepackage{xcolor}
\usepackage{array}
\usepackage{fontawesome5}

\usepackage[hidelinks]{hyperref}

\usepackage{algorithm}
\usepackage[noend]{algpseudocode}

\newcommand{\R}{\mathbb{R}}
\newcommand{\C}{\mathbb{C}}
\newcommand{\X}{\mathfrak{X}}

\newcommand{\del}{\partial}
\newcommand{\dd}{\mathrm{d}}
\newcommand{\ii}{\mathrm{i}}
\DeclareMathOperator{\im}{im}

\makeatletter
\g@addto@macro\bfseries{\boldmath}
\makeatother

\newtheorem{theorem}{Theorem}[section]
\newtheorem{corollary}[theorem]{Corollary}
\newtheorem{lemma}[theorem]{Lemma}
\newtheorem{proposition}[theorem]{Proposition}

\newtheorem{notation}{Notation}

\theoremstyle{definition}
\newtheorem{definition}[theorem]{Definition}
\newtheorem{example}[theorem]{Example}
\newtheorem{remark}[theorem]{Remark}

\addto\captionsUSenglish{}

\makeatletter
\newcommand{\MultiState}[1]{\parbox[t]{\dimexpr\linewidth-\ALG@thistlm\relax}{#1}}
\makeatother

\title{Geodesics on algebraic varieties}
\author[1]{Claire de Korte}
\author[2]{Yulia Mukhina} 
\author[3]{Paula Naomi Pilatus}

\affil[1]{Max Planck Institute for Mathematics in the Sciences, \href{mailto:claire.dekorte@mis.mpg.de}{claire.dekorte@mis.mpg.de}}

\affil[2]{ Department of Mathematical Sciences, University of Copenhagen, \href{mailto:yumu@math.ku.dk}{yumu@math.ku.dk}}

\affil[3]{Department of Mathematics, KU Leuven, \href{mailto:paulanaomi.pilatus@kuleuven.be}{paulanaomi.pilatus@kuleuven.be}}

\date{}

\begin{document}
\maketitle
\begin{abstract}
    This work initiates the study of geodesics from an applied algebraic perspective. We consider both affine and projective varieties, rephrasing the condition for a curve on these to be a geodesic in terms of the vanishing of the minors of an augmented Jacobian. We also provide analogous criteria for a curve to be a pregeodesic. 
    The Rosenfeld--Gr\"obner decompositions of the resulting systems of differential polynomials are studied, providing a proof of principle for the identification of geodesics on varieties through entirely algebraic methods. We complement this symbolic approach with a numerical algorithm, implemented in \texttt{Julia}, that computes geodesics between two given points, together with bounds on their distance.
\end{abstract}

\section{Introduction}

The study of geodesics dates back to the end of the $17$th century, when Johann Bernoulli asked for the shortest line between two points on a surface. An answer was supplied by Euler, in the form of a system of differential equations such lines must satisfy -- a set of \emph{geodesic equations} \cite{Euler1732}. Since then, geodesics have featured prominently in mathematical research, playing a key role in, for example, Gauss' theory of surfaces~\cite{Gauss1828}. Today, geodesics appear in a wide variety of applications, from optimal transport \cite{Mielke2013}, to protein dynamics~\cite{diepeveen2024riemannian}, to black hole physics \cite{SchwarzschildGeodesics}. Their broad utility continues to motivate research in mathematics and the sciences.

In practice the space on which geodesics are studied is often given by the zero set of a system of polynomials, that is an \emph{algebraic variety}. Despite this, most of the research on geodesics is formulated in the language of differential geometry, and does not leverage available algebraic tools. The study of \emph{metric algebraic geometry} \cite{metricalggeom} was initiated to bridge the gap between modern differential and algebraic geometry, translating insightful concepts from differential geometry into an algebraic language. In this paper, we expand on this and formulate necessary and sufficient conditions for a curve on an algebraic variety to be a geodesic in entirely (differential) algebraic language.

For instance, consider an affine variety $X=\mathcal{V}(F_1,\dots,F_k)\subset\mathbb{K}^n$, where $\mathbb{K}\in\{\R,\C\}$ and with metric induced by the Euclidean metric on $\mathbb{K}^n$. Let $c=\operatorname{codim}X$, Then a curve $\gamma\colon [a,b]\rightarrow X$ (which does not intersect the singular locus of~$X$) is a geodesic if and only if the $(c+1)\times(c+1)$-minors of the matrix  
\begin{equation}\label{eq: intro AJ}
A\mathcal{J}_{\mathbf{F}}(\gamma(t)):=\begin{pmatrix}\gamma''(t)\\ \vspace{-3mm}\\\overline{\mathcal{J}_{\mathbf{F}}(\gamma(t))}\end{pmatrix}, \end{equation} vanish for all $t\in[a,b]$. Here $\mathbf{F}=(F_1,\dots,F_k)$ and $\mathcal{J}_{\mathbf{F}}$ denotes the Jacobian of $\mathbf{F}$ and the bar refers to complex conjugation. We call the resulting system of equations the \emph{geodesic equation in implicit form}. 
    
We also derive a similar criterion for curves on projective varieties with metric induced by the Fubini--Study metric. For (affine or projective) varieties admitting a rational parametrization, we also state the geodesic equations in terms of this parametrization in classical form.

Our formulations bear striking resemblance to existing constructions in algebraic distance optimization, namely in the computation of \emph{Euclidean distance degrees} \cite{Draisma2015, friedman2026}. Given a (smooth) real affine variety $X=\mathcal{V}(\mathbf{F})\subset \R^n$, a generic point $\mathbf{u}\in \R^n$, the Euclidean distance degree is the number of critical points $\mathbf{x}\in X$ (seen as a complex variety) for the squared Euclidean distance $d_\mathbf{u}(\mathbf{x})=\sum_{i=1}^n(x_i-u_i)^2$. In~\cite{Draisma2015} it is shown that these critical points are given by the vanishing locus of the $(c+1)\times (c+1)$-minors of $$\begin{pmatrix}
    \mathbf{u}-\mathbf{x}\\
    \mathcal{J}_{\mathbf{F}}(\mathbf{x})
\end{pmatrix}.$$
There is a similar formulation of Euclidean distance degree for projective varieties~\cite{Draisma2015}.

Another feature of our implicit formulation of the geodesic equations of algebraic varieties is that these take the form of differential polynomials. Let us illustrate this using the example of the two-sphere $S^2=\mathcal{V}(x_1^2+x_2^2+x_3^2-1)\subset \R^3$.  For a curve~$\gamma(t)=(x_1(t),x_2(t),x_3(t))$ on $S^2$ the matrix $A\mathcal{J}_{\mathbf{F}}$ in~\eqref{eq: intro AJ} takes the shape $$\begin{pmatrix}
       x_1'' & x_2'' & x_3''\\
       2x_1 & 2x_2& 2x_3
\end{pmatrix}.$$
The curve $\gamma$ is a geodesic if and only if the $2\times 2$-minors of this matrix vanish, which yields the system
\begin{align*}
     x_1''x_2-x_1x_2''&=0,\\
x_1''x_3-x_1x_3''&=0,\\
      x_2''x_3-x_2x_3''&=0,
    \end{align*}
of polynomial differential equations. Similarly, when expressed in terms of a local rational parametrization, the geodesic equation in classical form takes the shape of differential polynomials after clearing denominators. This means that we can leverage the tools of differential algebra to study their properties and solutions. We note that geodesics on pseudo-Riemannian varieties were studied using differential algebra in \cite{Jaoui2019} from an abstract, model-theoretic perspective.

The differential-algebraic tool we use is the \emph{Rosenfeld--Gr\"obner} (RG) algorithm~\cite{BLOP1995}. In short, this algorithm splits the solutions of the system of polynomial differential equations into finitely many families, each described by a system of differential equations of a simple triangular form. Such decompositions give insight into the types of solutions to a set of differential equations, and, in some cases, can enable the construction of explicit solutions. As a proof of principle, we apply the RG algorithm to the geodesic equation of a number of different varieties, and thereby construct formulae for particular geodesics lying on them. These results are particularly stark if the identified geodesics are algebraic curves. In this case, the polynomials defining these curves appear directly in the output of the algorithm. 
The RG algorithm is thus a means by which one can identify algebraic geodesics without needing to solve differential equations. This perspective provides the foundation for formulating the conditions on a curve being a geodesic as an ideal membership problem.

Recall that geodesics are necessarily parametrized parallel to arc length, i.e.\ their velocity has constant norm. On the other hand many quantities such as the length of a curve are in fact reparametrization invariant. In many situations it is therefore enough to identify curves which are geodesics up to reparametrization, so-called \emph{pregeodesics}. Motivated by this, we provide generalizations of the criteria mentioned above to pregeodesics. For an affine variety $X=\mathcal{V}(\mathbf{F})\subset\mathbb{K}^n$ of codimension $c$, and a curve $\gamma\colon [a,b]\rightarrow X\setminus X_{\mathrm{sing}}$ this criterion is given by the vanishing of the $(c+2)\times (c+2)$-minors of the matrix obtained by appending $\gamma'(t)$ as an additional row to $A\mathcal{J}_{\mathbf{F}}(\gamma(t))$ as in~\eqref{eq: intro AJ} and the vanishing of the imaginary part of the scalar product of~$\gamma'$ and $\gamma''$. Again the resulting equations are differential polynomial and can therefore be studied using the RG algorithm.

In most applications the motivation to study (pre)geodesics stems from the fact that these are critical points for the length functional of paths with fixed endpoints. They can therefore be used to identify the shortest paths between a pair of points $\mathbf{p}$ and $\mathbf{q}$. The length of these shortest paths corresponds to the \emph{distance} between the points.

Computing the geodesic between two points $\mathbf{p}$ and $\mathbf{q}$ on a
Riemannian manifold, as well as the distance along it, is a classical problem. On smooth manifolds, one usually computes a geodesic from $\mathbf{p}$ to $\mathbf{q}$ by solving a
\emph{boundary value problem}. This is done numerically by shooting~\cite{Bryner2017,Sutti2024}, by finite differences~\cite{Maekawa1996}, or by
minimizing the energy of a discretized curve~\cite{KlassenSrivastavaMioJoshi2004,SrivastavaKlassen2016}. Most other numerical approaches use triangle meshes. Distances can then be computed efficiently by fast marching
or by heat flow~\cite{CraneWeischedelWardetzky2013,KimmelSethian1998}. For specific manifolds
such as spheres, Grassmannians and matrix
groups, there are also numerical implementations in software packages such as
Geomstats~\cite{Geomstats} and
Manopt~\cite{Manopt}. Complementing these packages we
provide software, written in \texttt{Julia}~\cite{Bezanson2017}, that computes
geodesics on algebraic varieties. 
As input, our implementation only needs the defining
equations of a variety or a rational
parametrization of it. 
The code
is available at
\begin{center}{\color{teal}\url{https://github.com/ymukhina/geodesics-on-algebraic-varieties}.}
\end{center}

This paper is structured as follows. In Section~\ref{sec: geodesics background} we introduce the notion of a (pre)geodesic on an affine or projective variety. In Sections~\ref{subsubsec: affine implicit geodesics} and~\ref{subsubsec: geodesics implicitly proj} we then state the geodesic equation in implicit form for affine and projective varieties, respectively. These are necessary and sufficient conditions for a curve on the variety to be a (pre)geodesic, phrased in terms of the defining polynomials of the variety. For varieties with a (local) rational parametrization we furthermore state the geodesic equations in terms of this parametrization in Section~\ref{subsubsec: geodesics in terms of parametrizations}. In Section~\ref{sec:precol}
we introduce tools from differential algebra and apply them in Section~\ref{subsec: special solutions} to study (pre)geodesics on certain varieties. 
More specifically, we apply the Rosenfeld--Gröbner algorithm to study special solutions to the (pre)geodesic equations in both implicit and parametric form. Finally, in Section~\ref{sec: numerical computation}, we present algorithms to compute the geodesic distance of a given pair of points on a variety alongside the geodesic realizing the distance. We test  the algorithm for several varieties using both the implicit and the parametric form of the geodesic equations.

\section{Geodesics on algebraic varieties}\label{sec: geodesics background}

In this section, we will introduce the notion of geodesics and discuss necessary and sufficient conditions for a curve on an algebraic variety to be a geodesic. In the following, we consider both affine varieties in $\mathbb{K}^n$ and projective varieties in $\mathbb{KP}^n$, where $\mathbb{K}=\R$ or $\mathbb{K}=\mathbb{C}$. 
\begin{notation} 
For $\mathbf{x}\in\mathbb{K}^n$, the Jacobian matrix of $\mathbf{F}=(F_1,\dots,F_k)$ at $\mathbf{x}=(x_1,\dots,x_n)$ is
\[
  \mathcal{J}_{\mathbf{F}}(\mathbf{x}) = \begin{pmatrix}
    \partial F_1/\partial x_1(\mathbf{x}) & \cdots & \partial F_1/\partial x_n(\mathbf{x})\\
    \vdots & \ddots & \vdots\\
    \partial F_k/\partial x_1(\mathbf{x}) & \cdots & \partial F_k/\partial x_n(\mathbf{x})
  \end{pmatrix}.
\]
For a complex matrix
$A = (a_{ij})$ we write $\overline{A} = (\overline{a_{ij}})$, where the bar denotes complex conjugation, and $A^{\dag}=\overline{A}^T$. 
\end{notation}
\begin{notation}
   For $\mathbb{K}=\R$, we denote the Euclidean inner product by
\begin{align*}
    \langle \mathbf{x},\mathbf{y}\rangle:=\sum_{i=1}^nx_iy_i,\quad \mathbf{x},\mathbf{y}\in\R^n.
\end{align*}
If $\mathbb{K}=\C$ we denote the Hermitian product by
\begin{align*}
    \langle \mathbf{x},\mathbf{y}\rangle:=\sum_{i=1}^nx_i\overline{y_i},\quad \mathbf{x},\mathbf{y}\in\mathbb{C}^n.
\end{align*}
Moreover, we denote the norm by
\begin{align*}
    \|\mathbf{x}\|:=\sqrt{\langle \mathbf{x},\mathbf{x}\rangle}, \quad \mathbf{x}\in\mathbb{K}^n.
\end{align*}
\end{notation}
In the following we will be focusing mostly on smooth algebraic varieties $X$, i.e.\ we will assume that $X=\mathcal{V}(\mathbf{F})$, where $\mathbf{F}=(F_1,\dots,F_k)$ are polynomials such that the Jacobian $\mathcal{J}_{\mathbf{F}}(\mathbf{x})$ has constant rank $c=\mathrm{codim} X$ for all $\mathbf{x}\in X$. However, we remark that our results can be applied in more general contexts to curves on singular varieties which do not meet the singular locus.

\subsection{Geodesics} 
Note that through the identification $T_{\mathbf{x}}\mathbb{K}^n\cong\mathbb{K}^n$, the inner product $\langle\cdot,\cdot\rangle$ defines a metric on the tangent space at every $\mathbf{x}\in\mathbb{K}^n$. 
We furthermore introduce a canonical metric on $T_{\mathbf{[x]}}\mathbb{KP}^n$ at every $[\mathbf{x}]\in\mathbb{KP}^n$ in the following way.

\begin{definition}    
\label{def: formula FS}
    The Fubini--Study metric at $[\mathbf{x}]\in\mathbb{KP}^n$ is defined by
$$g_{\mathrm{FS},[\mathbf{x}]}=\frac{\|\mathbf{x}\|^2\|\dd\mathbf{x}\|^2-(\mathbf{\bar{x}}\dd\mathbf{x})(\mathbf{x}\dd\mathbf{\bar{x}})}{\|\mathbf{x}\|^4},$$
where $\|\mathbf{x}\|^2=\sum_{i=0}^n x_i\bar{x}_i$, $\|\dd\mathbf{x}\|^2=\sum_{i=0}^n \dd x_i\dd\bar{x}_i$ and $\mathbf{\bar{x}}\dd\mathbf{x}=\sum_{i=0}^n\bar{x}_i\dd x_i$. 
\end{definition}

In the following, let
$X\subset \mathbb{K}^n$ or $X\subset\mathbb{KP}^n$ be a smooth affine or projective variety. 
We say that a continuos curve $\gamma\colon [a,b]\longrightarrow X$ is \emph{admissible} if there is a partition $(a_0,...,a_k)$ of $[a,b]$ such that for all $0\leq i<k$, the restriction $\gamma|_{[a_i,a_{i+1}]}$ is smooth and $\gamma'(t)\neq 0$ for all $t\in (a_i,a_{i+1})$. In this case we will say that $(a_0,..,a_k)$ is an \emph{admissible partition for} $\gamma$.

\begin{definition}
    Let $\gamma \colon [a,b]\longrightarrow X$ be an admissible curve and $(a_0,..,a_k)$ be an admissible partition of $\gamma$. The \emph{arc length} $L(\gamma)$ is defined as follows:
    \begin{itemize}
        \item If $X$ is affine, then \begin{equation*}
         L(\gamma):= \int_a^b \|\gamma'(t)\|\dd t:=\sum_{i=0}^{k-1} \int_{a_i}^{a_{i+1}} \|\gamma'(t)\|\dd t.   
        \end{equation*}
   \item If $X$ is projective, then
   \begin{equation*}
       L(\gamma):=\sum_{i=0}^{k-1} \int_{a_i}^{a_{i+1}} \sqrt{g_{\mathrm{FS}}\left(\gamma'(t),\gamma'(t)\right)}\dd t.\end{equation*} 
    \end{itemize}
\end{definition}

One can show that the arc length $L(\gamma)$ is invariant under reparametrizations of $\gamma$. Certain choices of parametrizations are particularly convenient to compute the arc length.

\begin{definition}
    We say that an admissible curve $\gamma$ on $X$ is \emph{parametrized parallel to arc length} if there is a constant $C>0$ such that whenever $\gamma'$ exists $\|\gamma'(t)\|=C$ for $X$ affine or $\sqrt{g_{\mathrm{FS}}\left(\gamma'(t),\gamma'(t)\right)}=C$ if $X$ is projective, respectively. If $C=1$, we say that $\gamma$ is \emph{parametrized by arc length}.
\end{definition}
Every admissible curve admits a reparametrization by arc length \cite[Prop 2.49]{Lee2018}. Note that if $\gamma\colon [a,b]\rightarrow X$ is parametrized by arc length then $L(\gamma)=b-a$. Moreover, if $\gamma\colon [0,1]\rightarrow X$ is parametrized parallel to arc length such that $\|\gamma'(t)\|=C$ or $\sqrt{g_{\mathrm{FS}}\left(\gamma'(t),\gamma'(t)\right)}=C$ whenever $\gamma'$ exists, then $L(\gamma)=C$. 
We will now state the definition of a \emph{geodesic}.
\begin{definition}
    An admissible curve $\gamma\colon [a,b]\rightarrow X$ is called a \emph{pregeodesic} if it is a critical point of the arc length functional $L(\gamma)$ (for variations of curves with fixed endpoints). It is called a \emph{geodesic} if it is a pregeodesic and parametrized parallel to arc length.
\end{definition}

Geodesics are necessarily smooth curves (cf.\ \cite[Thm 6.4]{Lee2018}). The geodesics on the (two-)sphere are great circles (cf.\ Example~\ref{ex: S2}). On the triaxial ellipsoid 
$$E_{a,b,c}:=\left\{\mathbf{x}\in\R^3:\tfrac{x_1^2}{a^2}+\tfrac{x_2^2}{b^2}+\tfrac{x_3^2}{c^2}=1\right\},\quad a,b,c\in\R\setminus\{0\},$$
geodesics are in general not closed. However, there exist closed geodesics, the simplest ones being the intersection with the three coordinate planes (cf.\ Example~\ref{ex: triaxial ellipsoid}).

Given an admissible curve $\gamma\colon [a,b]\longrightarrow X$, that curve being a geodesic is equivalent to it satisfying a system of second order ODEs. In Sections~\ref{subsubsec: affine implicit geodesics} and~\ref{subsubsec: geodesics implicitly proj} we will present this system of ODEs in terms of the polynomials defining $X$. Subsequently, in Section~\ref{subsubsec: geodesics in terms of parametrizations} we will see an alternative description in terms of a local (rational) parametrization of $X$.

\subsection{Geodesic equation in implicit form}
In the following, we will determine criteria for an admissible curve on an algebraic variety to be a geodesic in terms of its defining polynomials. We will consider the cases of affine and projective varieties separately in Sections~\ref{subsubsec: affine implicit geodesics} and~\ref{subsubsec: geodesics implicitly proj}. In both cases we show that the criterion can be phrased in terms of the vanishing of minors of certain augmented Jacobians. We provide similar necessary and sufficient criteria for an admissible curve to be a pregeodesic. 
\subsubsection{Affine varieties}
\label{subsubsec: affine implicit geodesics}

We will now state a criterion for a curve on an affine variety to be a geodesic in terms of the defining polynomials of the variety.

\begin{theorem}\label{thm:impl affine}
Let $X = \mathcal{V}(\mathbf{F})\subset \mathbb{K}^n$ be a smooth affine variety of codimension $c$, and let $\gamma\colon[a,b]\to\mathbb{K}^n$ be an admissible ($C^2$-)curve. Then $\gamma$ is a geodesic on $X$ if and only if for all $t\in[a,b]$, $\mathbf{F(}\gamma(t)) = 0$ and the matrix
 $$A\mathcal{J}_{\mathbf{F}}(\gamma(t)):=\begin{pmatrix}\gamma''(t)\\ \vspace{-3mm}\\\overline{\mathcal{J}_{\mathbf{F}}(\gamma(t))}\end{pmatrix}$$
 has vanishing $(c+1)\times(c+1)$-minors.
\end{theorem}
It is worth pointing out the striking similarities of Theorem~\ref{thm:impl affine} to the expressions appearing when computing the Euclidean distance degree of an algebraic variety~\cite[Eq.\ 2.1]{Draisma2015}. We call  $A\mathcal{J}_{\mathbf{F}}(\gamma(t))$ the \emph{augmented Jacobian}. 
\begin{proof}[Proof of Theorem~\ref{thm:impl affine}]
    A curve  $\gamma(t)=(x_1(t),...,x_n(t))$ is a geodesic on $X$  if and only if \cite[Prop 5.3]{Lee2018}
\[\mathbf{F}(x_1(t),...,x_n(t))=0\quad \text{and}\quad \gamma''(t)\perp T_{\gamma(t)}X.\]
 Since $T_{\gamma(t)}X= \ker \mathcal{J}_{\mathbf{F}}(\gamma(t))$ the condition becomes $\gamma''(t)\perp \ker \mathcal{J}_{\mathbf{F}}(\gamma(t))$ which is equivalent to the requirement that 
\[
\gamma''(t)\in\operatorname{span}\{\overline{\operatorname{grad}F_1(\gamma(t))},...,\overline{\operatorname{grad}F_k(\gamma(t))}\},
\]
where $\operatorname{grad}F_i=(\del F_i/\del x_1,\dots,\del F_i/\del x_n)$. This is equivalent to the vanishing of the $(c+1)\times (c+1)$-minors of $A\mathcal{J}_{\mathbf{F}}(\gamma(t))$ as defined in Theorem~\ref{thm:impl affine}.
\end{proof}
 $\mathbf{x}$ and $\mathbf{u}$. 

\begin{definition}
   We refer to the system of (differential) equations 
   $$ \mathbf{F}(\gamma)=0\quad \text{and}\quad A\mathcal{J}_{\mathbf{F}}(\gamma) \text{ has vanishing }(c+1)\times(c+1)\text{-minors}$$
   obtained in Theorem~\ref{thm:impl affine} as the \emph{geodesic equation (in implicit form)}.  
\end{definition}
Note that if $\mathbb{K}=\R$ the $(c+1)\times (c+1)$-minors of $A\mathcal{J}_{\mathbf{F}}(\gamma)$ are polynomials in the ring $\R[x_1,x''_1,\dots,x_n,x''_n]$. Similarly, if $\mathbb{K}=\C$, after splitting into real and imaginary parts, the $(c+1)\times (c+1)$-minors are polynomials in those and their second derivatives.

\begin{corollary}
    The geodesic equation in implicit form on $X$ is a system polynomial differential equations of order two, with polynomials of total degree less than or equal to $1-k+\sum_{i=1}^k\deg F_i$. If $\mathbb{K}=\R$, there are $\binom{n}{c+1}\cdot\binom{k+1}{c+1}$ equations from the condition on the minors and $k$ coming from the constraint $\mathbf{F}(\gamma)=0$. If $\mathbb{K}=\C$ there are $2\cdot\binom{n}{c+1}\cdot\binom{k+1}{c+1}$ equations from the condition on the minors and $2k$ coming from the constraint~$\mathbf{F}(\gamma)=0$. 
\end{corollary}

\begin{corollary}\label{cor: explicit}
Let $X = \mathcal{V}(\mathbf{F})\subset\mathbb{K}^n$ be a smooth affine variety of codimension $c$ with
$\mathbf{F} = (F_1,\dots,F_c)$. Then there are rational functions $R_i(\mathbf{x},\bar{\mathbf{x}},\mathbf{v})$, $i=1,\dots,n$, defined
for all $\mathbf{x}\in X$, $\mathbf{v}\in T_{\mathbf{x}}X$, such that an admissible ($C^2$-)curve $\gamma(t)=(x_1(t),\dots,x_n(t))$ in $X$ is a geodesic if and only if
\[
  \mathbf{x}'' = \mathbf{R}(\mathbf{x},\bar{\mathbf{x}},\mathbf{x}').
\]
Explicitly, the functions $R_i$ are given in terms of the Jacobian and the Hessian by
\begin{equation} \label{eq:formula-for-fun}
R_i(\mathbf{x},\bar{\mathbf{x}},\mathbf{v})=-\sum_{j,k}\mathcal{J}_{\mathbf{F}}(\mathbf{x})_{ij}^{\dag}\left(\mathcal{J}_{\mathbf{F}}(\mathbf{x})\mathcal{J}_{\mathbf{F}}(\mathbf{x})^{\dag}\right)^{-1}_{jk}(\mathbf{x}'^TH_{F_k}(\mathbf{x})\,\mathbf{x}').
\end{equation}
In particular, the
$R_i$ can be computed from $\mathbf{F}$ (and its first and second derivatives) by
linear algebra over the field of rational
functions in $\mathbf{x}$, $\bar{\mathbf{x}}$ and $\mathbf{v}$.
\end{corollary}
\begin{proof}
Suppose that $\gamma(t)=(x_1(t),\dots x_n(t))$ satisfies the geodesic equations in implicit form.
Then in particular the $(c+1)\times (c+1)$-minors of $A\mathcal{J}_{\mathbf{F}}(\gamma)$ vanish and hence there exist functions $\lambda_j\colon [a,b]\rightarrow \mathbb{K}$ such that
\begin{equation}\label{eq: cor gamma''=R}
    x_i''= \sum_{j=1}^c \lambda_j\,\overline{\tfrac{\del F_j}{\del x_i}(\mathbf{x})}.
\end{equation}
Since $F\circ \gamma=0$, it follows that
\begin{equation*}
    0=\tfrac{\dd^2}{\dd t^2}F_i\circ\gamma=\sum_{j}\tfrac{\del F_i}{\del x_j}(\mathbf{x})\,x_j''+\sum_k\tfrac{\del^2F_i}{\del x_j\del x_k}(\mathbf{x})\,x_j'x_k'.
\end{equation*}
Combining this with~\eqref{eq: cor gamma''=R}, we find that
\begin{equation*}
    -\sum_{j,k}\tfrac{\del^2F_i}{\del x_j\del x_k}(\mathbf{x})\,x_j'x_k'=\sum_{j,k}\tfrac{\del F_i}{\del x_j}(\mathbf{x})\overline{\tfrac{\del F_k}{\del x_j}(\mathbf{x})}\,\lambda_k=\sum_{k}\left(\mathcal{J}_{\mathbf{F}}(\mathbf{x})\mathcal{J}^{\dag}_{\mathbf{F}}(\mathbf{x})\right)_{ik}\lambda_k.
\end{equation*}
Since $\mathcal{J}_{\mathbf{F}}(\mathbf{x})$ has full rank for all $\mathbf{x}\in X$, the matrix $\mathcal{J}_{\mathbf{F}}(\mathbf{x})\mathcal{J}^{\dag}_{\mathbf{F}}(\mathbf{x})$ is invertible and we can solve for $\lambda_j$. Denoting by $H_{F_i}(\mathbf{x})$ the Hessian of $F_i$, we obtain
\begin{equation*}
    \lambda_i=-\sum_{k,j,\ell}\left(\mathcal{J}_{\mathbf{F}}(\mathbf{x})\mathcal{J}_{\mathbf{F}}(\mathbf{x})^{\dag}\right)^{-1}_{i\ell}H_{F_{\ell},jk}(\mathbf{x})\,x_j'x_k'= -\sum_{\ell}\left(\mathcal{J}_{\mathbf{F}}(\mathbf{x})\mathcal{J}_{\mathbf{F}}(\mathbf{x})^{\dag}\right)^{-1}_{i\ell}(\mathbf{x}'^TH_{F_{\ell}}(\mathbf{x})\,\mathbf{x}').
\end{equation*}
Plugging this back into~\eqref{eq: cor gamma''=R} yields
\begin{align*}
    x_i''&=-\sum_{j,k}\mathcal{J}_{\mathbf{F}}(\mathbf{x})_{ij}^{\dag}\left(\mathcal{J}_{\mathbf{F}}(\mathbf{x})\mathcal{J}_{\mathbf{F}}(\mathbf{x})^{\dag}\right)^{-1}_{jk}(\mathbf{x}'^TH_{F_k}(\mathbf{x})\,\mathbf{x}').
\end{align*}
Since the $F_i$ are polynomials, the expression on the right-hand side is a rational function in $\mathbf{x}$, $\bar{\mathbf{x}}$ and $\mathbf{x}'$. 
\end{proof}

\begin{remark}
For $\mathbb{K}=\R$, Corollary~\ref{cor: explicit} gives the classical form of the geodesic equations that are
usually derived with Lagrange multipliers, as the equations of motion of a free particle constrained to $X$ (see~\cite[Sec.~VII.1]{HairerLubichWanner2006} for constrained mechanical systems,
\cite[Sec.~10.2.2]{PatrikalakisMaekawa2002} for implicit surfaces and~\cite[Sec.~2.2.2]{EdelmanAriasSmith1998} for the
Stiefel manifold). We emphasize that Corollary~\ref{cor: explicit} requires that $X$ is defined by the vanishing of precisely $c=\operatorname{codim} X$ polynomials while Theorem~\ref{thm:impl affine} applies more generally.
\end{remark}

\begin{example}\label{ex: sphere}
    We consider the two-sphere $S^2\subset \R^3$, which is defined by the vanishing of
    $$ F(x_1,x_2,x_3)=x_1^2+x_2^2+x_3^2-1.$$
    Then the augmented Jacobian for a curve $\gamma(t)=(x_1(t),x_2(t),x_3(t))\in S^2\subset \R^3$ becomes $$A\mathcal{J}_F(\gamma(t))=\begin{pmatrix}
        x_1'' & x_2'' & x_3''\\
        2x_1 & 2x_2& 2x_3
    \end{pmatrix},$$
    where we left the $t$ dependence implicit. By Theorem~\ref{thm:impl affine} the curve $\gamma(t)$ is a geodesic if and only if it is a solution to the system of equations
    \begin{equation} \label{eq:for-the-example}
    \begin{aligned}
    x_1^2+x_2^2+x_3^2-1&=0,\\
    x_1''x_2-x_1x_2''&=0,\\
    x_1''x_3-x_1x_3''&=0,\\
    x_2''x_3-x_2x_3''&=0.
    \end{aligned}
    \end{equation}
    We will discuss in Section~\ref{sec: tools from differential algebra} how these equations can be studied using tools from differential algebra.
\end{example}

The following generalization of Theorem~\ref{thm:impl affine} gives necessary and sufficient conditions for an admissible curve to be a pregeodesic.
\begin{theorem}\label{thm:impl affine pregeodesic}
Let $X = \mathcal{V}(\mathbf{F})\subset \mathbb{K}^n$ be a smooth affine variety of codimension $c$, and let~$\gamma\colon[a,b]\to\mathbb{K}^n$ be an admissible
($C^2$-)curve. Then $\gamma$ is a pregeodesic on $X$ if and only if for all $t\in[a,b]$, $\mathbf{F}(\gamma(t)) = 0$, $\operatorname{Im} \langle \gamma''(t),\gamma'(t)\rangle=0$ and the matrix
 $$\begin{pmatrix}\gamma''(t)\\ \gamma'(t)\\\vspace{-3.5mm}\\\overline{\mathcal{J}_{\mathbf{F}}(\gamma(t))}\end{pmatrix} $$
 has vanishing $(c+2)\times(c+2)$-minors.
\end{theorem}
\begin{proof}
Suppose $\gamma\colon [a,b]\to \mathbb{K}^n$ is a reparametrization of a geodesic on $X$, i.e.\ there is a geodesic $\hat{\gamma}$ and a smooth function $\varphi\colon [a,b]\rightarrow \mathbb{R}$ such that $\gamma=\hat{\gamma}\circ\varphi$. Then $\gamma'(t)=\varphi'(t)\hat{\gamma}'(\varphi(t))$, \[\gamma''(t)=\varphi'(t)^2\, \hat{\gamma}''(\varphi(t))+\varphi''(t)\, \hat{\gamma}'(\varphi(t))\in \operatorname{span}_{\mathbb {K}}\{\gamma',\overline{\operatorname{grad}F_i}\},\]  
and \[\operatorname{Im}\langle \gamma'',\gamma'\rangle=\operatorname{Im}(\varphi'\varphi''\langle\hat{\gamma}'\circ\varphi,\hat{\gamma}'\circ\varphi\rangle)=0.\]

Conversely, suppose that $\gamma\colon [a,b]\to \mathbb{K}^n$ satisfies $F\circ \gamma=0$, $\operatorname{Im}\langle \gamma'',\gamma'\rangle=0$ and $\gamma''\in \operatorname{span}_{\mathbb {K}}\{\gamma',\overline{\operatorname{grad}F_i}\}$. Then by the same computation as above every reparametrization $\hat{\gamma}=\gamma\circ\tilde{\varphi}$ of $\gamma$ satisfies again $\hat{\gamma}''\in \operatorname{span}_{\mathbb {K}}\{\gamma',\overline{\operatorname{grad}F_i}\}$ and $\operatorname{Im}\langle \hat{\gamma}'',\hat{\gamma}'\rangle=0$. It remains to show that we can find a reparametrization such that~$\langle \hat{\gamma}'',\hat{\gamma}'\rangle=0$. To show this, define $f\colon [a,b]\rightarrow \R$ via
$f:=\langle\gamma',\gamma'\rangle$. Then we need to solve
\[0\overset{!}{=}\langle \hat{\gamma}'',\hat{\gamma} \rangle=\tfrac{1}{2}\tfrac{\dd}{\dd t}\langle \hat{\gamma}',\hat{\gamma}'\rangle=\tfrac{1}{2}\tfrac{\dd}{\dd t}\langle \tilde{\varphi}'\gamma'\circ\tilde{\varphi},\tilde{\varphi}'\gamma'\circ\tilde{\varphi}\rangle=\tfrac{1}{2}\tfrac{\dd}{\dd t}\tilde{\varphi}'^2\,f\circ\tilde{\varphi}=\tilde{\varphi}'\tilde{\varphi}''\,f\circ\tilde{\varphi}+\tfrac{1}{2}\tilde{\varphi}'^3\,f'\circ \tilde{\varphi}.\]

A solution $\tilde{\varphi}$ is given by
$\tilde{\varphi}^{-1}=C\int\sqrt{f}$
for some constant $C\in\R\setminus\{0\}$.
\end{proof}
\begin{remark} 
    Note that it follows from the second part of the proof of Theorem~\ref{thm:impl affine pregeodesic} that given a pregeodesic~$\gamma\colon [a,b]\rightarrow X$ a reparametrization $\varphi$ such that $\hat{\gamma}:=\gamma\circ \varphi$ is a geodesic is given by
    $$\varphi^{-1}(t)=\int_a^t \sqrt{\langle \gamma'(s),\gamma'(s)\rangle}\,\dd s.$$
    This yields a parametrization by arc length.
\end{remark}

\begin{example}\label{ex: pregeodesic sphere}
    We consider again the two-sphere $S^2\subset\mathbb{R}^3$ from Example \ref{ex: sphere}, except this time construct the pregeodesic equations. Via Theorem \ref{thm:impl affine pregeodesic}, these are
    \begin{equation}
        x_1^2+x_2^2+x_3^2-1=0,\quad\det\begin{pmatrix}x''_1&x_2''&x_3''\\x_1'&x_2'&x_3'\\2x_1&2x_2&2x_3\end{pmatrix}=0.
    \end{equation}
\end{example}

\subsubsection{Projective varieties}
\label{subsubsec: geodesics implicitly proj}
We will now suppose that $X\subset \mathbb{KP}^n$ is a smooth projective variety and state a theorem analogous to Theorem~\ref{thm:impl affine}. To this end, we need to introduce \emph{horizontal lifts}.
\begin{lemma}\label{lem: hor lift}
    Let $\gamma(t)$ be a curve in $\mathbb{KP}^n$. There is a curve $\tilde{\gamma}(t)$ on $S^N$ such that $\pi_{\mathbb{KP}^n}\circ\tilde{\gamma}=\gamma$ and such that $\langle\tilde{\gamma}(t),\tilde{\gamma}'(t)\rangle=0$. Moreover, $\tilde{\gamma}$ is unique up to a factor of $q\in\mathbb{K}$, $q\overline{q}=1$. 
\end{lemma}
\begin{proof}
    This is a special case of a more general result for principal bundles, see e.g.~\cite[Thm 5.8.2]{Hamilton2017}. Note that if $\mathbb{K}=\R$, the condition $\langle\tilde{\gamma}(t),\tilde{\gamma}'(t)\rangle=0$ is automatically satisfied.
\end{proof}

\begin{definition}
    The curve $\tilde{\gamma}$ as in Lemma~\ref{lem: hor lift} is called a \emph{horizontal lift} of $\gamma$.
\end{definition}

Since the proof of the following theorem requires some tools from differential geometry that we prefer not to introduce in the main part of this paper, we will discuss it in Appendix~\ref{appendix projective}.
\begin{theorem}\label{impl form proj : thm}
   Let $X = \mathcal{V}(\mathbf{F})\subset\mathbb{KP}^n$ be a smooth projective variety of codimension $c$, let $\gamma\colon[a,b]\to\mathbb{KP}^n$ be an admissible
($C^2$-)curve, and let $\tilde\gamma$ be a horizontal lift of $\gamma$. Then $\gamma$ is a geodesic of $X$ if and only if, for
all $t\in[a,b]$, $\mathbf{F}(\tilde\gamma(t)) = 0$ and all $(c+2)\times(c+2)$ minors of the matrix 
$$A_{\mathbb{KP}^n}\mathcal{J}_{\mathbf{F}}(\tilde{\gamma}(t)):=\begin{pmatrix}\tilde\gamma''(t)\\  \tilde\gamma(t)\\ \vspace{-3.5mm}\\ \overline{\mathcal{J}_{\mathbf{F}}(\tilde\gamma(t))}\end{pmatrix},$$ vanish. The conditions do not depend on the choice of horizontal lift.
\end{theorem} 
\begin{proof}
    This is a direct consequence of Proposition~\ref{appendix projective impl : prop} and Lemma~\ref{lem: appendix for rank of augm jac} in the appendix. Since horizontal lifts are unique up to a constant factor in $\mathbb{K}$, the criterion is independent of this choice. 
\end{proof}
Again we want to stress the similarities to the computation of the Euclidean distance degree of projective varieties~\cite[Eq.\ 2.7]{Draisma2015}. Moreover, as already mentioned above, the criterion also applies to curves outside the singular locus of a non-smooth projective variety.
\begin{remark}\label{rem: formulation proj impl}
    In order to identify all possible solutions to the geodesic equation, it will be useful to incorporate the condition that $\tilde{\gamma}(t)\in S^N$ by adding the equation $\|\tilde{\gamma}(t)\|^2=\sum_{i=0}^nx_i(t)\bar{x}_i(t)=1$ to the system. Moreover, we incorporate the condition of $\tilde{\gamma}$ being horizontal. In total the system becomes 
    \begin{align}\label{eq: rem on proj impl}
        \|\tilde{\gamma}\|^2=1,\quad \langle\tilde{\gamma},\tilde{\gamma}'\rangle=0, \quad \mathbf{F}(\tilde{\gamma})=0\quad\text{and}\quad A_{\mathbb{KP}^n}\mathcal{J}_{\mathbf{F}}(\tilde{\gamma}) \text{ has vanishing }(c+2)\times(c+2)\text{-minors}.
    \end{align}
    For $\mathbb{K}=\R$ the system~\eqref{eq: rem on proj impl} is the same as the one obtained when applying Theorem~\ref{thm:impl affine} to the affine variety in $\R^{n+1}$ given by $\mathcal{V}(\mathbf{F})\cap \mathcal{V}(\sum_i x_i^2 -1)$.
    For $\mathbb{K}=\C$ the system~\eqref{eq: rem on proj impl} is not algebraic due to the appearance of complex conjugation. It can be made algebraic by splitting into real and imaginary parts. 
\end{remark}

\begin{example}
    The Grassmannian $\mathrm{Gr}(2,4)$ can be seen as a codimension $1$ subvariety of $\mathbb{CP}^5$, defined by the Pl\"ucker relation
    \begin{equation*}
        [z_0:z_1:z_2:z_3:z_4:z_5],\quad z_0z_5-z_1z_4+z_3z_2=0.
    \end{equation*}
    Suppose that $\gamma(t)$ is a curve on the Grassmannian, with horizontal lift $\tilde{\gamma}(t)=(z_0(t),\ldots,z_5(t))$, seen as a curve in $\mathbb{C}^{6}$. Then $\tilde{\gamma}$ satisfies the equations
    \begin{equation*}
        z_0z_5-z_1z_4+z_3z_2=0, \quad z_0\bar{z}_0+\cdots+z_5\bar{z}_5=1,\quad z'_0\bar{z}_0+\cdots+z_5'\bar{z}_5=0,
    \end{equation*}
    and $\gamma$ is a geodesic if and only if the $(3\times3)$-minors of the matrix
    \begin{equation*}
        A_{\mathbb{CP}^5}\mathcal{J}_{\mathbf{F}}=\begin{pmatrix}z_0''&z_1''&z_2''&z_3''&z_4''&z_5''\\z_0&z_1&z_2&z_3&z_4&z_5\\\bar{z}_5&-\bar{z}_4&\bar{z}_3&\bar{z}_2&-\bar{z}_1&\bar{z}_0\end{pmatrix},
    \end{equation*}
    vanish. For this system to be encoded as differential polynomials, one must split it into real and imaginary parts. This doubles the number of equations in the system. 
\end{example}

Like in the case of affine varieties, Theorem~\ref{impl form proj : thm} can be generalized to give necessary and sufficient conditions for an admissible curve to be a pregeodesic.

\begin{theorem}\label{impl form proj pregeod: thm}
   Let $X = \mathcal{V}(\mathbf{F})\subset\mathbb{KP}^n$ be a smooth projective variety of codimension $c$, let $\gamma\colon[a,b]\to\mathbb{KP}^n$ be an admissible
($C^2$-)curve, and let $\tilde\gamma$ be a horizontal lift of $\gamma$. Then $\gamma$ is a pregeodesic of $X$ if and only if, for
all $t\in[a,b]$, $\mathbf{F}(\tilde\gamma(t)) = 0$, $\operatorname{Im}\langle \tilde{\gamma}''(t),\tilde{\gamma}'(t)\rangle=0$ and all $(c+3)\times(c+3)$ minors of the matrix 
$$\begin{pmatrix}\tilde\gamma''(t)\\ \tilde{\gamma}'(t)\\\tilde\gamma(t)\\ \vspace{-3.5mm}\\ \overline{\mathcal{J}_{\mathbf{F}}(\tilde\gamma(t))}\end{pmatrix},$$ vanish. The conditions do not depend on the choice of horizontal lift.
\end{theorem} 
\begin{proof}
    The proof works largely parallel to that of Theorem~\ref{thm:impl affine pregeodesic}. To see that  $$\operatorname{rank}\begin{pmatrix} \tilde{\gamma}'(t)\\\tilde\gamma(t)\\ \vspace{-3.5mm}\\ \overline{\mathcal{J}_{\mathbf{F}}(\tilde\gamma(t))}\end{pmatrix}=c+2,$$
    use that $\tilde{\gamma}'$ and $\tilde{\gamma}$ are orthogonal to each other and both vectors are orthogonal to $\operatorname{span}\{\overline{\operatorname{grad}F_i}\}$ (cf.\ Lemma~\ref{lem: appendix for rank of augm jac}). Moreover, the conditions are invariant under different choices of horizontal lifts, since these are related by a constant factor in $\mathbb{K}$.
\end{proof}

\subsection{Geodesics in terms of rational parametrizations}\label{subsubsec: geodesics in terms of parametrizations}

Let $X$ be a smooth affine or projective variety. In the following we consider the case where $X$ admits a local (real) rational parametrization of the algebraic variety $X$, i.e.\ open subsets $U\subset \mathbb{R}^m$, $V\subset X$ and a rational invertible map $\varphi\colon U\rightarrow V$. Here, $m=\dim_{\R}X$. 
\begin{definition}
    Let $X$ be a smooth algebraic variety and $\varphi\colon U\rightarrow V$ a local rational parametrization of $X$. The \emph{induced metric} $h\colon U\rightarrow \mathbb{R}^{m\times m}$ is defined as follows.
    \begin{itemize}
        \item If $X$ is affine: $h_{ij}(\mathbf{s}):=\operatorname{Re}\left\langle \frac{\del\mathbf{f}}{\del s_i}(\mathbf{s}),\frac{\del\mathbf{f}}{\del s_j}(\mathbf{s})\right\rangle$, where $\varphi(\mathbf{s})=(f_1(\mathbf{s}),\dots,f_n(\mathbf{s}))$;
        \item If $X$ is projective : $h_{ij}(\mathbf{s})=\operatorname{Re}g_{\mathrm{FS}}\left(\frac{\del\mathbf{f}}{\del s_i}(\mathbf{s}),\frac{\del\mathbf{f}}{\del s_j}(\mathbf{s})\right)$, where $\varphi(\mathbf{s})=[f_0(\mathbf{s}):f_1(\mathbf{s}):\dots:f_n(\mathbf{s})]$.     \end{itemize}
\end{definition}
The induced metric has the interpretation of restricting the metric of the ambient space to the tangent space of $X$ at $\varphi(\mathbf{s})$, $\mathbf{s}\in U$. In the following lemma, we state explicit formulas for $h(\mathbf{s})$.
\begin{lemma}Let $X$ be a smooth algebraic variety and $\varphi\colon U\rightarrow V$ a local rational parametrization of $X$.\begin{itemize}
    \item If $X$ is affine and $\varphi(\mathbf{s})=(f_1(\mathbf{s}),\dots,f_n(\mathbf{s}))$, then $h(\mathbf{s})=\tfrac{1}{2}\left(\mathcal{J}_{\mathbf{f}}(\mathbf{s})^T\cdot \overline{\mathcal{J}_{\mathbf{f}}(\mathbf{s})}+\mathcal{J}_{\mathbf{f}}(\mathbf{s})^\dag\cdot \mathcal{J}_{\mathbf{f}}(\mathbf{s})\right)$.
    \item If $X$ is projective and $\varphi(\mathbf{s})=[f_0(\mathbf{s}):f_1(\mathbf{s}):\dots:f_n(\mathbf{s})]$, then 
    \begin{equation}\label{eq: Fubini-Study induced metric real projective}h(\mathbf{s})=\tfrac{1}{2}\left(\|\mathbf{f}(\mathbf{s})\|^{-2}\, \mathcal{J}_{\mathbf{f}}(\mathbf{s})^T\overline{\mathcal{J}_{\mathbf{f}}(\mathbf{s})}-\|\mathbf{f}(\mathbf{s})\|^{-4}\,\left(\mathcal{J}_{\mathbf{f}}(\mathbf{s})^T\overline{\mathbf{f}(\mathbf{s})}\right)\otimes\left(\mathcal{J}_{\mathbf{f}}(\mathbf{s})^T\overline{\mathbf{f}(\mathbf{s})}\right)^{\dag}\quad + \text{c.c.}\right),\end{equation}
    where c.c.\ stands for complex conjugate.
\end{itemize}
    
\end{lemma}
\begin{proof} The affine case is clear. Let now $X\subset \mathbb{KP}^n$ be a smooth projective variety.
    Recall from Definition~\ref{def: formula FS} that the Fubini--Study metric at $[\mathbf{x}]\in\mathbb{KP}^n$ is given by 
$$g_{\mathrm{FS}}([\mathbf{x}]
)=\frac{\|\mathbf{x}\|^2\|\dd\mathbf{x}\|^2-(\mathbf{\bar{x}}\dd\mathbf{x})(\mathbf{x}\dd\mathbf{\bar{x}})}{\|\mathbf{x}\|^4}.$$
For the components of the induced metric we compute
\begin{align*}
    h_{ij}(\mathbf{s})&=\operatorname{Re}g_{\mathrm{FS}}\left(\tfrac{\del\mathbf{f}}{\del s_i}(\mathbf{s}),\tfrac{\del\mathbf{f}}{\del s_j}(\mathbf{s})\right)\\
    &=\operatorname{Re}\left(\|\mathbf{f}(\mathbf{s})\|^{-2}\sum_k \tfrac{\del f_k}{\del s_i}(\mathbf{s}) \overline{\tfrac{\del f_k}{\del s_j}(\mathbf{s})}- \|\mathbf{f}(\mathbf{s})\|^{-4}\left(\sum_k \overline{f_k(\mathbf{s})}\tfrac{\del f_k}{\del s_i}(\mathbf{s})\right)\left(\sum_l f_l (\mathbf{s})\overline{\tfrac{\del f_l}{\del s_j}(\mathbf{s})}\right)\right)\\
    &=\tfrac{1}{2}\left(\|\mathbf{f}(\mathbf{s})\|^{-2}\, \left(\mathcal{J}_{\mathbf{f}}(\mathbf{s})^T\overline{\mathcal{J}_{\mathbf{f}}(\mathbf{s})}\right)_{ij}-\|\mathbf{f}(\mathbf{s})\|^{-4}\,\left(\mathcal{J}_{\mathbf{f}}(\mathbf{s})^T\overline{\mathbf{f}(\mathbf{s})}\right)_i\left(\mathcal{J}_{\mathbf{f}}(\mathbf{s})^\dag\mathbf{f}(\mathbf{s}\right)_j\quad+\text{c.c.}\right).
    \end{align*}
\end{proof}

Every curve $\gamma\colon [a,b]\rightarrow V$ that lies in $V\subset X$ can be expressed in terms of this local parametrization
\[ \mathbf{s}(t):=\varphi^{-1}\circ \gamma(t):=(s_1(t),\dots,s_m(t)).\]
Given such a parametrization, the condition that $\gamma$ is a geodesic can be phrased in terms of a coupled system of ODEs for the functions $s_i(t)$ as follows.

We collect the derivatives of the components of $h(\mathbf{s})$ in the following $(m\times m\times m)$-tensor \[h_{ij,k}(\mathbf{s}):=\tfrac{\del h_{ij}}{\del s_k}(\mathbf{s}),\]
and use it to define the \emph{Christoffel symbols}
\begin{align}\label{eq: christoffel}C_{jkl}(\mathbf{s})=\frac{1}{2}\sum_{i=1}^mh^{-1}_{ji}(\mathbf{s})(h_{ik,l}(\mathbf{s})+h_{il,k}(\mathbf{s})-h_{kl,i}(\mathbf{s})).\end{align}

\begin{theorem}[See e.g.\ Eq 4.16 in \cite{Lee2018}] \label{geod eq param form: thm}
   An admissible ($C^2$-)curve $\mathbf{s}(t)=(s_1(t),\dots s_m(t))$ in $U$ is mapped under $\varphi$ to a geodesic in $V$ if and only if the functions $s_i(t)$ solve the following system of ODEs:
\begin{equation}
\label{eq: geodesic eq in local coordinates}
    s''_j(t)+\sum_{k,l=1}^m C_{jkl}(\mathbf{s}(t))\, s'_k(t)s'_l(t)=0.
\end{equation} 
\end{theorem}
\begin{definition}
 We refer to \eqref{eq: geodesic eq in local coordinates} as \emph{geodesic equation (in parametric form)}.   
\end{definition}
\begin{remark}
   We observe that (after clearing denominators) the left hand side of equation \eqref{eq: geodesic eq in local coordinates} are polynomials in the ring $\mathbb{R}[s_1,s_1',s_1'',\ldots,s_m,s_m',s_m'']$.  
\end{remark}
 
For the sake of completeness we will include a direct proof that for geodesics on affine varieties contained in the image of a local parametrization, the system in parametric form in Theorem~\ref{geod eq param form: thm} is equivalent to the geodesic equation in implicit form stated in Theorem~\ref{thm:impl affine}. 
\begin{lemma}
    Let $X=\mathcal{V}(F_1,...,F_k)\subset \mathbb{K}^n$ be a smooth affine variety and $\varphi\colon (U\subset \mathbb{R}^m)\rightarrow (V\subset X)$ a local rational parametrization of $X$. If a curve $\mathbf{s}$ in $U$ satisfies the geodesic equation in parametric form then $\gamma:=\varphi\circ \mathbf{s}$ satisfies the geodesic equation in implicit form. Conversely, for every curve $\gamma\subset V$ satisfying the geodesic equation in implicit form, $\mathbf{s}:=\varphi^{-1}\circ\gamma$ satisfies the geodesic equation in parametric form.
\end{lemma}
\begin{proof}
    We only show this for the case $\mathbb{K}=\R$, since the complex affine varieties can be treated as real affine varieties of twice the dimension. Recall that the implicit geodesic equations are equivalent to the requirement that $\gamma''(t)\perp \ker \mathcal{J}_{\mathbf{F}}(\gamma(t))$ for every $t$, where $\mathbf{F}=(F_1,\dots,F_k)$. Since $\varphi$ is a parametrization of $V(\mathbf{F})$, we have $\mathbf{F}\circ \varphi=0$ and hence~$\mathcal{J}_{\mathbf{F}}\cdot \mathcal{J}_{\varphi}=0$. By the assumptions that $X$ is smooth and $\varphi$ is a local parametrization, we know that $\mathcal{J}_{\mathbf{F}}$ is of rank $c:=n-m$ and $\mathcal{J}_{\varphi}$ is of rank $m$ so we have $\im \mathcal{J}_{\varphi}= \ker \mathcal{J}_{\mathbf{F}}\mid_{\operatorname{im}(\varphi)}$.
    Therefore, proving that $\gamma''(t)$ is perpendicular to $\ker \mathcal{J}_{\mathbf{F}}$ is equivalent to showing that 
    \[\left\langle \gamma''(t),\tfrac{\del \mathbf{f}}{\del s_j}(\mathbf{s}(t))\right\rangle=0,\quad \text{for}\quad j=1,\dots,m.\]
    To this end, we 
    compute
    \begin{equation}\label{eq: impl equiv to param}
    \begin{aligned}
        \left\langle \gamma'',\tfrac{\del \mathbf{f}}{\del s_j}\right\rangle&=\sum_{i,k,l}\tfrac{\del f_i}{\del s_j}\tfrac{\del^2 f_i}{\del s_k\del s_l}s'_k s'_l+\sum_{i,k} \tfrac{\del f_i}{\del s_j}\tfrac{\del f_i}{\del s_k}s''_k\\
        &=\sum_{i,k,l}\tfrac{\del f_i}{\del s_j}\tfrac{\del^2 f_i}{\del s_k\del s_l}s'_k s'_l+\sum_{k} h_{jk}s''_k.\end{aligned}
    \end{equation}
    Substituting \eqref{eq: geodesic eq in local coordinates}, the second term becomes
    \begin{align*}
        \sum_{k} h_{jk}s''_k&=-\sum_{i,k,l}h_{jk}C_{kil}s'_i s'_l\\
        &=-\tfrac{1}{2}\sum_{i,l}(h_{ji,l}+h_{jl,i}-h_{il,j})s'_i s'_l\\
        &=-\sum_{i,k,l}\tfrac{\del f_k}{\del s_j}\tfrac{\del^2f_k}{\del s_i\del s_l}s'_i s'_l.
    \end{align*}
    This shows that $\langle \gamma'',\del \mathbf{f}/\del s_j\rangle=0$ for all $j$ and hence $\gamma$ satisfies the implicit geodesic equation. Conversely, if $\langle \gamma'',\del \mathbf{f}/\del s_j\rangle=0$ by~\eqref{eq: impl equiv to param} it follows that
    \begin{align*}
     \sum_{k} h_{jk}\left(s''_k+\sum_{i,l}C_{kil}s_i's_l'\right)=0.   
    \end{align*}
    Since $h$ is invertible, this implies the geodesic equation in parametric form.
\end{proof}

\begin{example}\label{ex: two-sphere parametric}
    By inverting the stereographic projection, one finds the following rational parametrization of the two-sphere
    \begin{align*}
        f_1(s_1,s_2)=\frac{2s_1}{s_1^2+s_2^2+1},\quad f_2(s_1,s_2)=\frac{2s_2}{s_1^2+s_2^2+1},\quad
        f_3(s_1,s_2)=\frac{s_1^2+s_2^2-1}{s_1^2+s_2^2+1}.
    \end{align*}
    Following the procedure described above, we use this to construct the induced metric and Christoffel symbols. Using the formula \eqref{eq: geodesic eq in local coordinates}, we find the parametric geodesic equations to be
    \begin{align*}
        -4s_2s_1's_2'+2s_1(-(s_1')^2+(s_2')^2)+(1+s_1^2+s_2^2)s_1''&=0,\\-4s_1s_1's_2'+2s_2((s_1')^2-(s_2')^2)+(1+s_1^2+s_2^2)s_2''&=0.
    \end{align*}
    We will discuss in Section~\ref{sec: tools from differential algebra} how these equations can be studied using differential algebra.
\end{example}

As for the implicit equations, there is also a generalization of the geodesic equations in parametric form to pregeodesics.
\begin{theorem}\label{pregeod eq param form: thm}
   An admissible ($C^2$-)curve $\mathbf{s}(t)=(s_1(t),\dots s_m(t))$ in $U$ is mapped under $\varphi$ to a pregeodesic in $V$ if and only if there is some smooth function $f\colon [a,b]\to \R$ such that the functions $s_i(t)$ solve the following system of ODEs:
\begin{equation}
\label{eq: pregeodesic eq in local coordinates}
    s''_j(t)+\sum_{k,l=1}^m C_{jkl}(\mathbf{s}(t))\, s'_k(t)s'_l(t)=f(t)\, s'_j(t).
\end{equation} 
\end{theorem}

\section{Geodesics from differential algebra}\label{sec: tools from differential algebra}

In Section~\ref{sec: geodesics background} we demonstrated that, with respect to our choices of metric, systems of geodesic and pregeodesic equations of algebraic varieties can be written as differential polynomial equations. In this section, we exploit this fact, and leverage tools from differential algebra to analyze these systems and their solutions.

\subsection{Differential algebra preliminaries} \label{sec:precol}

We start by recalling the relevant background from
differential algebra, following
Ritt~\cite{Ritt1950} and
Kolchin~\cite{Kolchin1973}. We refer
to~\cite{boulier2023intro, BLOP1995,
BLOP2009, Hubert2003} for the description of
the Rosenfeld–Gröbner algorithm and the proofs of its
properties.

\begin{definition} A \emph{differential ring} $(R, \;')$ is a commutative ring with a derivation $': R \rightarrow R$, that is, a map such that, for all $a,b \in R$, $(a+b)' = a' +b'$ and $(ab)' = a'b + ab'$. A \emph{differential field} is a differential ring which is a field. An ideal $I \subset R $ is called a \emph{differential ideal} if $a' \in I$ for every $a \in I$.
\end{definition}

\begin{notation}
Let $x$ be an element of a differential ring. 
We write $x^{(\infty)} := \{x,x',x'',x^{(3)},\ldots\}.$
\end{notation}

\begin{definition} 
    Let $R$ be a differential ring, and consider the ring of
polynomials in infinitely many variables
\[
R[x^{(\infty)}] := R[x, x', x'', x^{(3)},\ldots].
\]
Extend the derivation from $R$ to this ring by $(x^{(j)})':= x^{(j+1)}$ \cite[§ 9, Prop. 4]{bourbaki1950elements}. The resulting differential ring is called the \emph{ring of differential polynomials in $x$ over $R$}.
The ring of differential polynomials in several variables is defined by iterating this construction.
\end{definition}

\begin{example}
    The geodesic equations for the sphere~\eqref{eq:for-the-example}
    can be thought of as a system of differential polynomials in the ring $\mathbb{R}[x_1^{(\infty)}, x_2^{(\infty)}, x_3^{(\infty)}]$.
\end{example}

\begin{definition}
    A \emph{ranking} is a total order $\prec$ on
the set $x_1^{(\infty)}\cup\cdots\cup
x_n^{(\infty)}$ of all derivatives
such
that, for all $1\le i,j\le n$ and all
$k,l,m \in \mathbb{Z}_{\ge0}$,
    \begin{itemize}
    \setlength\itemsep{-0.3em}
        \item $x_i^{(k)}\prec x_i^{(k+1)}$,
        \item $x_i^{(k)}\prec x_j^{(l)}\implies
        x_i^{(k+m)}\prec x_j^{(l+m)}$.
    \end{itemize}
\end{definition}

In what follows, we will mainly work with the
following two rankings, namely the elimination
and the orderly ranking.

\begin{definition}
    The \emph{elimination ranking}, denoted $[x_1,\ldots,x_n]$, is the
    ranking satisfying, for all $1 \le i,j\le n$
    and all $k,l\in \mathbb{Z}_{\ge0}$,
    \[
      x_i^{(k)}\succ x_j^{(l)}\quad\text{if }
      i<j,\text{ or if }i=j\text{ and }k>l.
    \]
    The \emph{orderly ranking}, denoted
$[[x_1,\ldots,x_n]]$, is the
    ranking satisfying, for all $1\le i,j\le n$
    and all $k,l\in \mathbb{Z}_{\ge0}$,
    \[
      x_i^{(k)}\succ x_j^{(l)}\quad\text{if }
      k>l,\text{ or if }k=l\text{ and }i<j.
    \]
\end{definition}

\begin{notation}
Following the notation of~\cite{DifferentialAlgebra}, by
$[[x_1,\ldots,x_i],x_{i+1},\ldots,x_n]$ we
mean the ranking that compares two derivatives
of $x_1,\ldots,x_i$ using the orderly ranking,
and all other pairs of derivatives using the
elimination ranking. 
\end{notation}

The central tool from differential algebra
that we will use is the so-called
\emph{Rosenfeld--Gr\"obner} (RG)
\emph{algorithm}~\cite{BLOP1995}.
Based on characteristic
sets~\cite{Ritt1950}, it works as follows.
Let $K$ be a differential field of characteristic zero, and $\Sigma=\{g_1,\ldots,g_s\}$ be a finite
set of differential polynomials in
$K[x_1^{(\infty)},\ldots,
x_n^{(\infty)}]$.  We
write $\Sigma=0$ for the system of
differential eqxuations $g_1=\cdots=g_s=0$. Let $I$ be the radical differential ideal generated by
$\Sigma$. Given a ranking, the RG algorithm
computes a decomposition
\begin{equation*}
I=J_1\cap\cdots\cap J_r,
\end{equation*}
into radical differential ideals $J_i$.
Each $J_i$ is generated by a certain finite
set $\Omega_i$ of differential polynomials, up to
a saturation. The sets $\Omega_i$
depend on the choice of ranking. For details,
we refer to~\cite{Hubert2003}. Up to this saturation, the solution set of
$\Sigma=0$ is the union of the solution sets
of the systems $\Omega_1=0,\ldots,\Omega_r=0$
(see~\cite[Sec.~3.9]{boulier2023intro}). We call $J_i$ \emph{generic} if its solutions depend on as many constants as the general solution of $\Sigma=0$, and \emph{non-generic} otherwise. For geodesics on a variety $X$, the general solution depends on $2\cdot\dim X$ constants: a starting point and an initial velocity. 
The sets $\Omega_i$ of
generic components often contain many large
differential polynomials. In contrast,
non-generic components represent `special'
solutions, whose differential equations can look
particularly simple in the chosen coordinates.

The
Rosenfeld–Gröbner algorithm has been implemented in the BLAD library~\cite{BLAD}, which underlies
the \texttt{DifferentialAlgebra} package of the Maple standard library and is also accessible from
SageMath and Python\footnote{\url{https://pypi.org/project/DifferentialAlgebra/}}.

\begin{example}
\label{ex: rosenfeld grobner sphere}
Recall the system of geodesic equations
arising from the implicit formulation of the
two-sphere,
\begin{equation*}
  \Sigma=\{x_1^2+x_2^2+x_3^2-1,\,\,
  x_1''x_2-x_1x_2'',\,\,
  x_1''x_3-x_1x_3'',\,\,
  x_2''x_3-x_2x_3''\}.
\end{equation*}
Applying the RG algorithm with ranking
$[[x_1,x_2,x_3]]$ to this system yields a
decomposition into three components
$J_1,J_2$ and $J_3$, described by the sets
    \begin{align*}
    \Omega_1 &= \{x_1^2 + x_2^2 + x_3^2 - 1,\,\,x_2^2 x_3 (x_3')^2 + x_2^2x_3''- 2x_2 x_3^2 x_2' x_3'+x_3^3(x_2')^2+ x_3^2 x_3''-x_3(x_2')^2- x_3(x_3')^2- x_3'',\\
    &\quad\quad\quad x_2^3 (x_3')^2- 2x_2^2 x_3 x_2' x_3'+ x_2^2 x_2'' + x_2x_3^2 (x_2')^2- x_2 (x_2')^2- x_2 (x_3')^2 + x_3^2 x_2'' - x_2''\},\\
    \Omega_2&=\{x_2^2+x_3^2-1,\,\,x_1,\,\, x_3^2x_3''-x_3(x_3')^2-x_3''\},\\
    \Omega_3&=\{x_3^2-1,x_2,x_1\}.
    \end{align*}
The component $J_1$ is generic, while $J_2$
and $J_3$ are non-generic. The system of
equations $\Omega_2=0$ encodes a `special'
solution to $\Sigma=0$, the circle
$x_2^2+x_3^2=1$ in the $x_1=0$ plane. Here the
differential equation
$x_3^2x_3''-x_3(x_3')^2-x_3''=0$ encodes the
parametrization of the circle that is
required for it to be a geodesic. It is
solved by $x_3(t)=\sin(c_1t+c_2)$. This implies that $x_2(t)=\pm\cos(c_1t+c_2)$, and so the curve $\gamma=(0,x_2,x_3)$ has constant speed $\|\gamma'(t)\|=|c_1|$. In particular, this tells us that the curve is indeed parametrised proportional to arc length, with start point specified by $c_2$.
The system of equations $\Omega_3=0$
corresponds to the constant solutions
$(0,0,1)$ and $(0,0,-1)$.
\end{example}
Example~\ref{ex: rosenfeld grobner sphere}
already demonstrates how the RG algorithm can
be useful for finding geodesics on algebraic
surfaces, particularly when these geodesics
are algebraic curves and we are using the
implicit formulation of the geodesic
equations. Further applications of the RG
algorithm to geodesics can be found in
Section~\ref{subsec: special solutions}.

\subsection{Special solutions via Rosenfeld--Gr\"obner}\label{subsec: special solutions}

Described in Section \ref{sec:precol}, the Rosenfeld--Gr\"obner algorithm is a powerful tool for understanding the types of solutions of a system of differential equations. Example \ref{ex: rosenfeld grobner sphere} already shows how the RG algorithm can be useful for determining equations for geodesics. In this section, we further demonstrate this, applying the RG algorithm to find specific geodesics on a number of different varieties. In each case, we compute geodesics on affine varieties with metric induced from the Euclidean metric on $\mathbb{R}^n$. These examples provide a proof of principle for the identification of geodesics using purely differential-algebraic manipulations. 
All decompositions were computed with~\texttt{DifferentialAlgebra}~5.1~%
\cite{DifferentialAlgebra}. The code that
reproduces these computations, together with
the raw output of all decompositions, can be found in the \texttt{differential-algebra} folder of the GitHub repository
\begin{center}
{\color{teal}\url{https://github.com/ymukhina/geodesics-on-algebraic-varieties}}
\end{center}

We start by revisiting the sphere from Example \ref{ex: rosenfeld grobner sphere}.

\begin{figure}[h!]
    \centering
    \begin{subfigure}[b]{0.24\textwidth}
        \centering
        \includegraphics[
            trim={2cm 2cm 2.5cm 0cm},
            clip,
            width=\linewidth
        ]{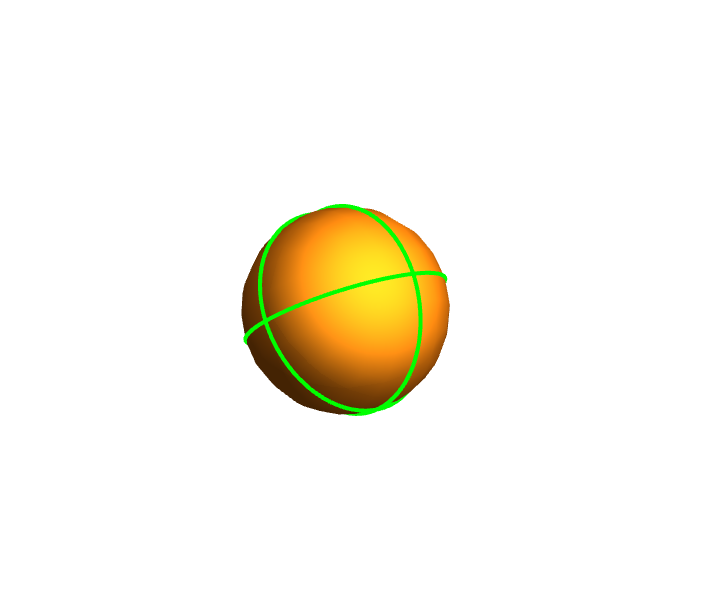}
        \caption{Two-sphere.}
        \label{fig: sphere geodesics}
    \end{subfigure}\hfill
    \begin{subfigure}[b]{0.24\textwidth}
        \centering
        \includegraphics[
            trim={3cm 3.3cm 2cm 2cm},
            clip,
            width=\linewidth
        ]{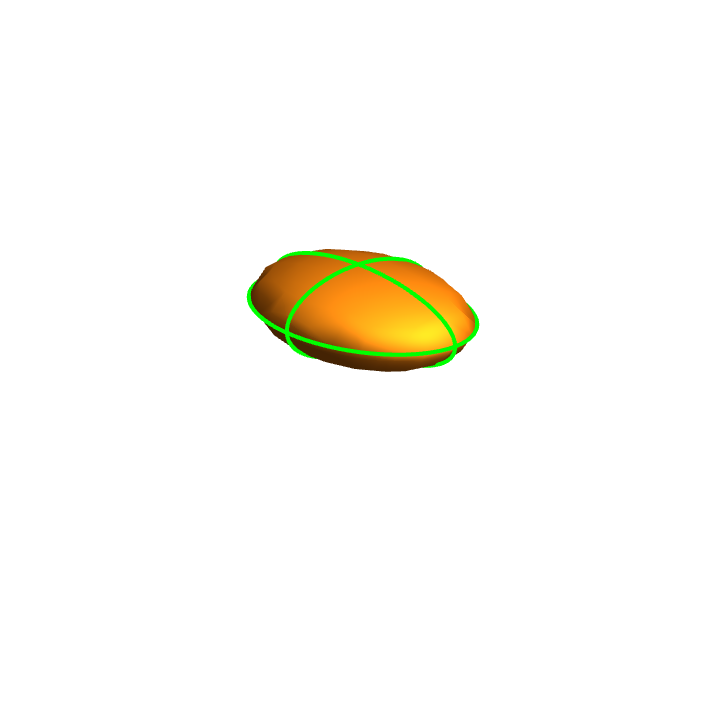}
        \caption{Triaxial ellipsoid.}
        \label{fig: triaxial ellipsoid geodesics}
    \end{subfigure}\hfill
    \begin{subfigure}[b]{0.24\textwidth}
        \centering
        \includegraphics[
            width=\linewidth
        ]{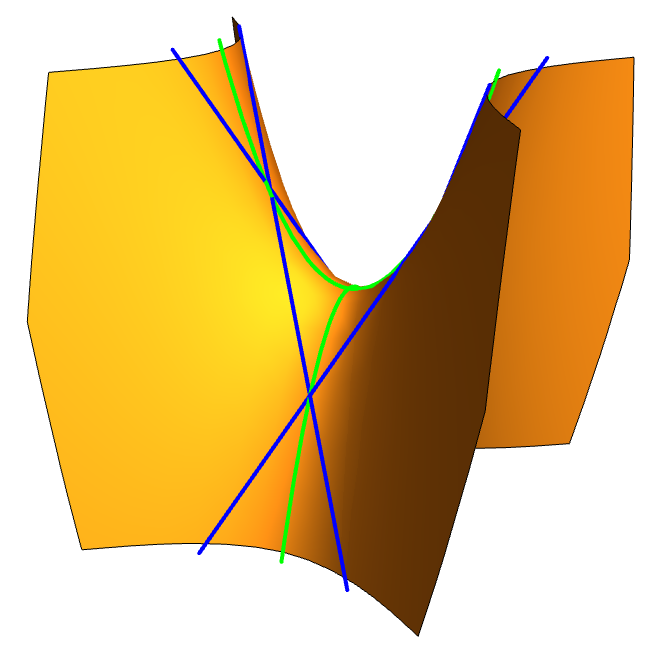}
        \caption{Hyperbolic paraboloid.}
        \label{fig: hyperbolic paraboloid geodesics}
    \end{subfigure}\hfill
    \begin{subfigure}[b]{0.24\textwidth}
        \centering
        \includegraphics[
            width=\linewidth
        ]{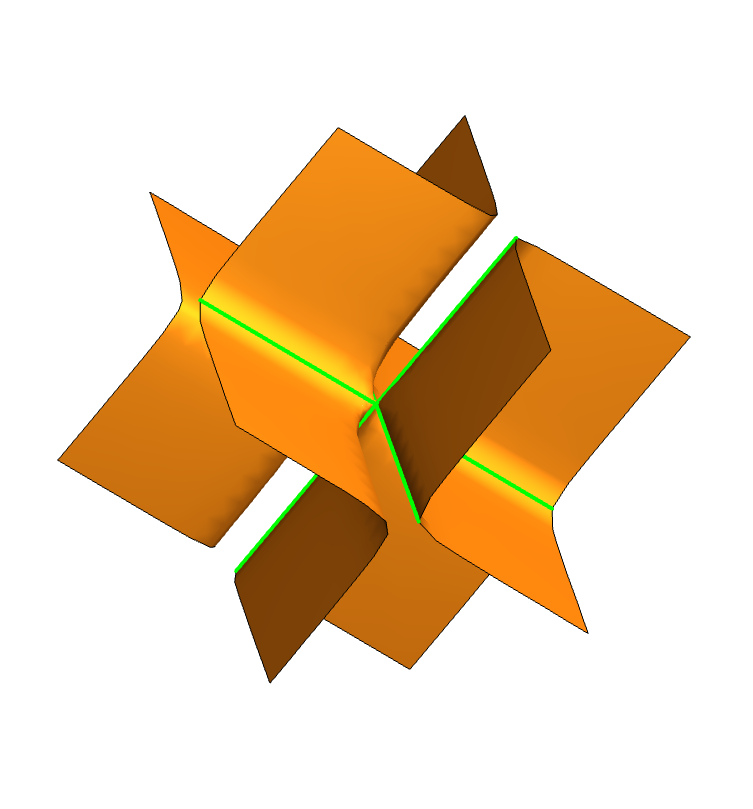}
        \caption{The Cayley cubic.}
        \label{fig: cayley cubic line}
    \end{subfigure}

    \caption{Examples of geodesics (green or blue) on surfaces in $\mathbb{R}^3$ (orange), all identified using the RG algorithm. 
    Those in blue belong to the family of lines defined in equation \eqref{eqn: lines on hyperbolic paraboloid}.}
    \label{fig: RG geodesic plots}
\end{figure}

\begin{example}[Sphere]\label{ex: S2}
    As described in Example \ref{ex: rosenfeld grobner sphere}, with ranking $[[x_1,x_2,x_3]]$, the RG algorithm allowed us to identify the circle $x_2^2+x_3^2=1$ in the $x_1=0$ plane as a geodesic. Since the system of geodesic equations is invariant under any permutation of the variable labels, that decomposition also tells us that the unit circles in the $x_2=0$ and $x_3=0$ planes are geodesics. This can also be seen without permuting any of the indices by applying the RG algorithm using the ranking $[x_2,x_1,x_3]$. One finds a decomposition consisting of $6$ components. Two of these components contain multiple large differential polynomials, two encode sets of points and two are `special'. The `special' components are specified by the differential polynomials
    \begin{align}\label{eqn: special sphere geodesic components}
        \Omega_1=\{x_2,x_1^2+x_3^2-1,x_3^2x_3''-x_3(x'_3)^2-x_3''\},\quad
        \Omega_2=\{x_3,x_1^2+x_2^2-1,x_1^2x_1''-x_1(x_1')^2-x_1''\}.
    \end{align}
    These correspond to the circles of unit radius in the $x_2=0$ and $x_3=0$ planes respectively, parametrised in the same manner as the circle from Example \ref{ex: rosenfeld grobner sphere}. One could instead apply the RG algorithm to the pregeodesic equations from Example \ref{ex: pregeodesic sphere}. Using the ranking $[x_1,x_2,x_3]$, one finds a decomposition into $6$ components, one generic, three corresponding to points, and two `special' components, specified by
    \begin{equation}
        \Omega_1=\{x_2^2+x_3^2-1,x_1\},\quad \Omega_2=\{x_1^2+x_2^2-1,x_3\},
    \end{equation}
    which correspond to circles in the $x_1=0$ and $x_3=0$ planes respectively. By symmetry, using the ranking $[x_1,x_3,x_2]$, one additionally finds the circle in the $x_2=0$ plane. As expected, unlike the components in \eqref{eqn: special sphere geodesic components}, these don't contain any information about how these circles should be parametrised to get geodesics. The circles that the RG algorithm was able to identify from the geodesic (or pregeodesic) equations are plotted in Figure \ref{fig: sphere geodesics}.
\end{example}

This basic example demonstrates that different geodesics can be identified by different choices of ranking. We now move on to consider the classical example, first considered by Jacobi \cite{Jacobi1839}, of geodesics on ellipsoids.

\begin{example}[Triaxial Ellipsoid]\label{ex: triaxial ellipsoid}
    Consider the following triaxial ellipsoid, seen as a subvariety of $\mathbb{R}^3$
    \begin{equation*}
        \frac{x_1^2}{9}+\frac{x_2^2}{4}+x_3^2-1=0.
    \end{equation*}
    One can compute the geodesic equations for this variety using the implicit formulation of Theorem \ref{thm:impl affine}. Applying the RG algorithm to the resulting system, using the ranking $[[x_1,x_2,x_3]]$, we find a decomposition into $6$ components. One component is generic, two correspond to individual points, two have empty solution sets (if we instead worked over $\mathbb{C}$, these would also just correspond to points), and one is `special'. The `special' component is specified by the differential polynomials
    \begin{equation*}
        \Omega=\{x_1,x_2^2+4x_3^2-4,3x_3^4x_3''-2x_3^2x_3''-4x_3(x_3')^2-x_3''\}.
    \end{equation*}
    The solution to $\Omega=0$ is an ellipse in the $x_1=0$ plane, parametrised in such a way that $x_3$'s dependence on $t$ takes the form of the inverse function of an Elliptic integral of the second kind. The ranking $[[x_2,x_1,x_3]]$ yields an analogous decomposition, except this time the `special' component is specified by
    \begin{equation*}
        \Omega=\{x_1^2+9x_3^2-9,x_2,8x_3^4x_3''-7x_3^2x_3''-9x_3(x_3')^2-x_3''\}.
    \end{equation*}
   Likewise for the ranking $[[x_3,x_1,x_2]]$, with `special' component specified by
   \begin{equation*}
       \Omega=\{4x_1^2+9x_2^2-36,x_3,5x_2^4x_2''-4x_2^2x_2''-36x_2(x_2')^2-64x_2''\}.
   \end{equation*}
    These correspond to ellipses in the $x_2=0$ and $x_3=0$ planes respectively. As for the ellipse in the $x_1=0$ plane, the differential polynomials in these sets tell us that for these ellipses to be geodesics, they must be parametrised in terms of the inverse function of an elliptic integral of the second kind. As we did with the sphere, we could equally apply the RG algorithm to the pregeodesic equations for this variety, derived using Theorem \ref{thm:impl affine pregeodesic}. Using the ranking $[x_1,x_2,x_3]$, we find a decomposition into $6$ components, one generic, three consisting of individual points, and two `special', specified by
    \begin{equation*}
        \Omega_1=\{x_2^2+4x_3^2-4,x_1\},\quad\Omega_2=\{4x_1^2+9x_2^2-36,x_3\},
    \end{equation*}
    which correspond to the same ellipses in the $x_1=0$ and $x_3=0$ planes that we found earlier. Using the ranking $[x_1,x_3,x_2]$ we find the ellipse in the $x_2=0$ plane too. 
    The ellipses identified by the RG algorithm from the (pre)geodesic equations are plotted in Figure \ref{fig: triaxial ellipsoid geodesics}.
\end{example}

Another family of simple surfaces is the hyperbolic paraboloids. As we will see, the RG algorithm can also identify geodesics on these.

\begin{example}[Hyperbolic Paraboloid]
    Consider the following hyperbolic paraboloid, seen as an affine subvariety of $\mathbb{R}^3$
    \begin{equation*}
        x_1^2-x_2^2-2x_3=0.
    \end{equation*}
    Computing the geodesic equations corresponding to this variety in the implicit formulation of Theorem \ref{thm:impl affine}, and then applying the RG algorithm with ranking $[x_3,x_1,x_2]$, we find a decomposition into $5$ components, one generic, one consisting of individual points, and three `special'. The `special' components are specified~by
    \begin{equation*}
    \begin{split}
        &\Omega_1=\{x_2,x_1^2x_1''+x_1(x_1')^2+x_1'',-x_1^2+2x_3\},\\
        &\Omega_2=\{x_2^2x_2''+x_2(x_2')^2+x_2'',x_1,x_2^2+2x_3\},\\&\Omega_3=\{x_2'',(x_1')^2-(x_2')^2,-x_1^2+x_2^2+2x_3\},
    \end{split}
    \end{equation*}
    The solutions to $\Omega_1=0$ and $\Omega_2=0$ are parabolas in the $x_2=0$ and $x_1=0$ planes respectively, parametrised such that the $t$ dependence of $x_1$ or $x_2$ is given in terms of inverse hyperbolic functions. The general solution to $\Omega_3=0$ is the family of lines
    \begin{equation}\label{eqn: lines on hyperbolic paraboloid}
        x_1(t)=\pm at+c,\quad x_2(t)=at+b,\quad x_3(t)=(\pm c-b)at+\frac{1}{2}(c^2-b^2).
    \end{equation}
    For each choice of sign, and each choice of real $a,b$ and $c$, we have a geodesic. Using the same ranking on the pregeodesic equations, we find a decomposition into just $3$ components, one generic, one consisting of individual points, and one `special'. The special component is specified by
    \begin{equation*}
        \Omega=\{x_2,-x_1^2+2x_3\},
    \end{equation*}
    and simply corresponds to the parabola in the $x_2=0$ plane that we identified above. The parabola in the $x_1=0$ plane can similarly be found from the pregeodesic equations using the ranking $[x_3,x_2,x_1]$, and the lines defined by $x_1^2-x_2^2=0$ in the $x_3=0$ plane from the ranking $[[x_1,x_2],x_3]$. 
    The parabolas and a subset of the lines the RG algorithm identified as geodesics are plotted in Figure \ref{fig: hyperbolic paraboloid geodesics}. 
\end{example} 

The Cayley cubic is a projective variety in $\mathbb{RP}^3$ defined by the equation
\begin{equation*}
        x_1x_2x_3+x_1x_2x_0+x_2x_3x_0+x_1x_3x_0=0.
    \end{equation*}
It is singular, and is the only cubic surface with four nodes. Whilst a generic cubic surface will contain $27$ lines, the Cayley cubic contains just $9$ \cite{hunt1996geometry}. Instead of thinking of the Cayley cubic as a projective variety, we dehomogenise by setting $x_0=1$. As we will see in the example below, the RG algorithm identifies three lines on the resulting surface as geodesics.

\begin{example}[The Cayley Cubic]
    Now consider the Cayley cubic, seen as an affine subvariety in $\mathbb{R}^3$, cut out by
    \begin{equation}\label{eqn: Cayley in x0=1 chart}
        x_1x_2x_3+x_1x_2+x_2x_3+x_1x_3=0.
    \end{equation}
     Again employing the implicit formulation of Theorem \ref{thm:impl affine}, we construct an associated system of geodesic equations. Applying the RG algorithm to this with ranking $[[x_1,x_2,x_3]]$, we find a decomposition into $6$ components. One is generic, four consist of individual points and one is `special'. The `special' component is specified by the differential polynomials
    \begin{equation}
        \Omega=\{x_3,x_2,x_1''\}.
    \end{equation}
    The solution to $\Omega=0$ consists of a line along the $x_1$-axis, parametrised such that $x_1$ has linear dependence on $t$. Like with the sphere example, since equation \eqref{eqn: Cayley in x0=1 chart} is invariant under permutations of the variable labels, this tells us that the lines along the $x_2$ and $x_3$ axes are also geodesics (away from the origin). This can also be seen explicitly from the output of the RG algorithm, applying it to the geodesic equations using the rankings $[[x_2,x_1,x_3]]$ and $[[x_3,x_1,x_2]]$ respectively.
    These lines are plotted in Figure \ref{fig: cayley cubic line}. Unlike in the previous examples, regardless of the choice of ranking, the RG algorithm did not finish within our limits (cf.\ Remark~\ref{rem: run time}) when applied to the system of pregeodesic equations for this surface.
    \end{example}

 The previous examples were all considering systems of differential equations constructed using the implicit formulations of Theorems \ref{thm:impl affine} and \ref{thm:impl affine pregeodesic}, where the RG algorithm enabled us to find geodesics as algebraic curves. The following example demonstrates that the RG algorithm can also be useful in the parametric formulation.

\begin{example}[The Special Orthogonal Group]
    Consider the following rational parametrization of the special orthogonal group $\mathrm{SO}(3)$
    \begin{equation}
        \begin{pmatrix}f_1&f_2&f_3\\f_4&f_5&f_6\\f_7&f_8&f_9\end{pmatrix}=\frac{1}{1+s_1^2+s_2^2+s_3^2}\begin{pmatrix}  1-s_1^2-s_2^2+s_3^2&-2(s_1+s_2s_3)&2(-s_2+s_1s_3)\\2(s_1-s_2s_3)&1-s_1^2+s_2^2-s_3^2&-2(s_1s_2+s_3)\\2(s_2+s_1s_3)&2(-s_1s_2+s_3)&1+s_1^2-s_2^2-s_3^2
        \end{pmatrix}.
    \end{equation}
    This parametrization does not cover all of $\mathrm{SO}(3)$. Applying the process explained in Section \ref{subsubsec: geodesics in terms of parametrizations}, we construct the geodesic equations arising from this choice of parametrization, arriving at the system of ODEs
    \begin{align*}
    -2s_1'(s_1s_1'+s_2s_2'+s_3s_3')+(1+s_1^2+s_2^2+s_3^2)s_1''&=0,\\
    -2s_2'(s_1s_1'+s_2s_2'+s_3s_3')+(1+s_1^2+s_2^2+s_3^2)s_2''&=0,\\
    -2s_3'(s_1s_1'+s_2s'_2+s_3s'_3)+(1+s_1^2+s_2^2+s_3^2)s_3''&=0.
    \end{align*}
    Then applying the RG algorithm, using the ranking $[[s_1,s_2],s_3]$, we find a decomposition consisting of $13$ components. Eight of these components contain multiple large differential polynomials, three have empty solution sets over $\mathbb{R}$, and two are `special'. The two `special' components are specified by the differential polynomials
    \begin{align*}
        \Omega_1&=\{s_3^2s_3''-2s_3(s_3')^2+s_3'',s_1^2+s_2^2,s_3^2s_2''-2s_3s_2's_3'+s_2''\},\\
        \Omega_2&=\{s_1,s_2,s_3^2s_3''-2s_3(s_3')^2+s_3''\}.
    \end{align*}
    Since the $s_i$ are valued in $\mathbb{R}$, each of these components encodes the same solution space. This solution space corresponds to a straight line along the $s_3$-axis, parametrised such that $s_3$'s dependence on $t$ is given by $s_3(t)=\tan(c_1t+c_2)$. As in the case of the sphere, $c_2$ specifies the start point of the geodesic, and $c_1$ the norm of its velocity. Substituting this solution back into the parametrization, we find
    \begin{equation}
        \begin{pmatrix}
            1&0&0\\0&\cos(2(c_1t+c_2))&-\sin(2(c_1t+c_2))\\0&\sin(2(c_1t+c_2))&\cos(2(c_1t+c_2))
        \end{pmatrix},
    \end{equation}
    and so the RG algorithm has allowed us to identify rotations around the `$x$-axis' with constant angular velocity as geodesics.
\end{example}

\begin{remark}\label{rem: run time}
The running time and memory usage of the
\verb|DifferentialAlgebra| implementation of
the RG algorithm depend strongly on the
ranking. In the examples above, the
computation finished within a few seconds for
some rankings, while for others it did not
finish within our limits of $1$ hour and
$10$\,GB of memory per ranking.
\end{remark}

\section{Computing geodesics and geodesic distances}\label{sec: numerical computation}

Throughout this section, we work with a smooth connected real affine variety $X = \mathcal{V}(\mathbf{F}) = \mathcal{V}(F_1,\dots,F_c)\subset\mathbb{R}^n$ of codimension $c$. 
We note that in practice smoothness is only
needed along the computed geodesic on the
variety. Moreover, the instead of restricting to connected varieties, we can restrict to a particular connected component of a given variety. Note furthermore that complex affine varieties in
$\mathbb{C}^n$ and real projective varieties
in $\mathbb{RP}^n$ can be reduced to the
setting of real affine varieties in $\mathbb{R}^{2n}$ and
$\mathbb{R}^{n+1}$, respectively (see Remark~\ref{rem: formulation proj impl} for the case of real projective varieties).

The numerical methods employed in this section
are standard: we solve initial value problems
using a Runge--Kutta method and boundary value
problems by single shooting. Our aim is to
illustrate how the geodesic equations of the form as in
in Corollary~\ref{cor: explicit} can be used to
compute geodesics and geodesic distances, that
is, the lengths of shortest curves on $X$
between two given points. 
We therefore describe the general computational procedures for completeness, ensuring consistency between the theoretical framework of this work and its numerical implementation. Our implementation complements existing software packages such as Manopt~\cite{Manopt}, Geomstats~\cite{Geomstats}, and Manifolds.jl~\cite{ManifoldsJl}, by providing tools specifically designed for computing geodesics on algebraic varieties.

\subsection{Algorithm}\label{sec: distance}

In this section, we present an algorithm to
compute the \emph{geodesic distance}, which
we now define.
\begin{definition}
\label{def: geodesic distance}
    Given a pair of points $\mathbf{p}, \mathbf{q}\in X$, the \emph{geodesic distance} $\dd_X(\mathbf{p},\mathbf{q})$ between them is defined as the infimum over the arc lengths of admissible curves on $X$ connecting $\mathbf{p}$ and $\mathbf{q}$.
    An admissible curve is \emph{minimizing} if its length is the geodesic distance between
its endpoints.
\end{definition}

A minimizing curve is a critical point of
the arc length functional and hence, after
reparametrization, a geodesic. We therefore
look for geodesics from $\mathbf{p}$ to
$\mathbf{q}$. Algorithm~\ref{alg:geodesic}
computes such a geodesic, represented by a
tuple of points on it, together with its
length. Up to numerical errors, this length is an upper bound for
$\dd_X(\mathbf{p},\mathbf{q})$, with equality
if the geodesic is minimizing.
The algorithm combines classical numerical
methods: a Runge--Kutta method~\cite{HNW} to
solve the geodesic equations with a given
initial point and initial velocity, Newton's
method to move the computed points back
onto $X$~%
\cite[Sec.~IV.4]{HairerLubichWanner2006},
and single
shooting~\cite[Ch.~2]{Keller1968},
\cite[Sec.~7.3]{StoerBulirsch2002} with the
Gauss--Newton method to correct the initial
velocity until the endpoint is $\mathbf{q}$.

\begin{algorithm}[ht]
\caption{A geodesic from $\mathbf{p}$ to
$\mathbf{q}$ and its length}\label{alg:geodesic}

\begin{algorithmic}[1]
  \Require
   \begin{tabular}[t]{@{}l@{}}
    \textbullet~the polynomials
    $\mathbf{F} = (F_1,\dots,F_c)
    \in\mathbb{R}[x_1,\dots,x_n]^c$;\\
    \textbullet~points $\mathbf{p},\mathbf{q}
    \in\mathcal{V}(\mathbf{F})$;\\
    \textbullet~a number of steps
    $K\in\mathbb{Z}_{>0}$;\\
    \textbullet~a tolerance $\delta>0$ for the
    Gauss--Newton method;\\
    \textbullet~a tolerance $\varepsilon>0$;\\
    \textbullet~a maximal number $N$ of shots.
  \end{tabular}
 
   \Ensure points $\mathbf{P}_0,\dots,\mathbf{P}_K$
  on a geodesic from $\mathbf{p}$ to $\mathbf{q}$
  on $\mathcal{V}(\mathbf{F})$ and its length, or
  \texttt{fail}
  \State \label{line:equations} Compute $\mathbf{R} = (R_1,\dots,R_n)\in\mathbb{R}(\mathbf{x},\mathbf{x}')^n$ using~\eqref{eq:formula-for-fun} (see Corollary~\ref{cor: explicit})
  \State \label{line:guess} $\mathbf{v}\gets$ the orthogonal projection of $\mathbf{q}-\mathbf{p}$ onto $T_{\mathbf{p}}\mathcal{V}(\mathbf{F})$
  \State $t_k\gets k/K$ for $k$ from $0$ to $K$
  \Repeat \label{line:loop}
    \State $\mathbf{P}_0\gets\mathbf{p}$
    \State $\mathbf{w}_0\gets\mathbf{v}$
    \For{$k$ \textbf{from} $0$ \textbf{to} $K-1$}
      \State \label{line:shoot}\MultiState{$\mathbf{x}_{\mathbf{v}}\gets$ the solution of
        \[
          \mathbf{x}'' = \mathbf{R}(\mathbf{x},\mathbf{x}'),\quad \mathbf{x}(t_k) = \mathbf{P}_k,\quad
          \mathbf{x}'(t_k) = \mathbf{w}_k
        \]
        on $[t_k,t_{k+1}]$ (Runge--Kutta method)}
      \If{$\|\mathbf{F}(\mathbf{x}_{\mathbf{v}}(t_{k+1}))\|>\varepsilon$} \label{line:shootcheck}
        \State \Return \texttt{fail}
      \EndIf
      
       \State \label{line:shootproject} $\mathbf{P}_{k+1}\gets$ the point of $\mathcal{V}(\mathbf{F})$ found by Newton's
        method for $\mathbf{F} = 0$ started at $\mathbf{x}_{\mathbf{v}}(t_{k+1})$
      \State \label{line:shootprojectv} $\mathbf{w}_{k+1}\gets$ the orthogonal projection of
        $\mathbf{x}_{\mathbf{v}}'(t_{k+1})$ onto $T_{\mathbf{P}_{k+1}}\mathcal{V}(\mathbf{F})$
    \EndFor
    \If{$\|\mathbf{P}_K-\mathbf{q}\|\ge\delta$}
      \State \label{line:newton} Correct $\mathbf{v}\in T_{\mathbf{p}}\mathcal{V}(\mathbf{F})$ by a Gauss--Newton step
        for $\mathbf{P}_K = \mathbf{q}$~\cite[Chapter~4]{Deuflhard2004}
    \EndIf
  
  \Until{$\|\mathbf{P}_K-\mathbf{q}\|<\delta$
    or $N$ shots have been made}

  \If{$\|\mathbf{P}_K-\mathbf{q}\|\geq\delta$}
    \State \Return \texttt{fail}
  \EndIf

    \State \label{line:return} \Return $(\mathbf{P}_0,\dots,\mathbf{P}_K)$
    and $\|\mathbf{v}\|$
\end{algorithmic}
\end{algorithm}

\begin{remark}\label{rem: variants}
Algorithm~\ref{alg:geodesic} also works if $X$
is given by a parametrization
$\varphi\colon U\to X$. Let $\mathbf{p} = \varphi(\mathbf{s}_\mathbf{p})$ and
$\mathbf{q} = \varphi(\mathbf{s}_\mathbf{q})$. By
Theorem~\ref{geod eq param form: thm}, a curve
$\mathbf{s}(t)$ in $U$ is mapped by $\varphi$ to
a geodesic if and only if it solves the
geodesic equations in parametric form. So we
run Algorithm~\ref{alg:geodesic} in $U$, with
$\mathbf{s}_\mathbf{p}$ and $\mathbf{s}_\mathbf{q}$ in place of
$\mathbf{p}$ and $\mathbf{q}$, and with the following 
changes.
\begin{itemize}
  \item Line~\ref{line:equations}: we compute
    the geodesic equations in parametric form
    instead of computing $\mathbf{R}$.

  \item Line~\ref{line:return}: we return the
    points $\varphi(\mathbf{s}(t_k))$ and the
    length $\|(\varphi\circ\mathbf{s})'(0)\|$
    instead of $\|\mathbf{v}\|$.
\end{itemize}

\end{remark}

To make the steps of Algorithm~\ref{alg:geodesic} more transparent, we illustrate it in the following example.

\begin{example}\label{ex: sphere distance}
Let $\mathbf{F} = x_1^2+x_2^2+x_3^2-1$, and let $\mathbf{p} = (1,0,0)$,
$\mathbf{q} = (0,0,1)$. Thus, $X = \mathcal{V}(\mathbf{F}) \subset\mathbb{R}^3$ is the unit sphere. The shortest curve from $\mathbf{p}$ to $\mathbf{q}$ is a quarter of a great circle, so we expect
$\dd_X(\mathbf{p},\mathbf{q}) = \pi/2$.

In this example we let $K = 3$, $\varepsilon = 10^{-9}$, $\delta = 10^{-11}$ and $N = 20$. Then Algorithm~\ref{alg:geodesic} runs as follows:
\begin{enumerate}
    \item By Corollary~\ref{cor: explicit}, the geodesic equations $\mathbf{x}'' = \mathbf{R}(\mathbf{x},\mathbf{x}')$ are given by
        \[
          \left\{
          \begin{aligned}
            x_1'' &= R_1(\mathbf{x}, \mathbf{x}') := -\frac{x_1'^2+x_2'^2+x_3'^2}{x_1^2+x_2^2+x_3^2}\,x_1,\\[2pt]
            x_2'' &= R_2(\mathbf{x}, \mathbf{x}') := -\frac{x_1'^2+x_2'^2+x_3'^2}{x_1^2+x_2^2+x_3^2}\,x_2,\\[2pt]
            x_3'' &= R_3(\mathbf{x}, \mathbf{x}') := -\frac{x_1'^2+x_2'^2+x_3'^2}{x_1^2+x_2^2+x_3^2}\,x_3.
          \end{aligned}
          \right.
        \]
  \item The initial guess $\mathbf{v} = (0,0,1)$.
  \item  $t_0 = 0$, $t_1 = 1/3$, $t_2 = 2/3$, $t_3 = 1$.

 \item Starting from $\mathbf{P}_0 = \mathbf{p}$, lines~\ref{line:shoot}--\ref{line:shootprojectv} give the points
          \[
            \mathbf{P}_1\approx(0.945,\,0,\,0.327),\qquad \mathbf{P}_2\approx(0.786,\,0,\,0.618),\qquad
            \mathbf{P}_3\approx(0.540,\,0,\,0.841).
          \]
          on the great circle through $\mathbf{p}$ and $\mathbf{q}$ (see more details in Figure~\ref{fig: sphere shooting}(a)). We note that $\mathbf{P}_3$ misses $\mathbf{q}$, indeed
          $\|\mathbf{P}_3-\mathbf{q}\| = 0.56\ge\delta$. Moreover, the points
   $\mathbf{x}_{\mathbf{v}}(t_1)$,
   $\mathbf{x}_{\mathbf{v}}(t_2)$ and
   $\mathbf{x}_{\mathbf{v}}(t_3)$, appearing in
   the intermediate computation from the
   Runge--Kutta method, lie on the sphere or inside it, at distance less than
   $10^{-15}$. Newton's method in
   line~\ref{line:shootproject} moves each of
   them back onto the sphere, along the ray
   through the origin.

  \item Thus, we correct $\mathbf{v}$ as in line~\ref{line:newton} of Algorithm~\ref{alg:geodesic} and shoot again:
        \begin{center}
        \begin{tabular}{cll}
          \hline
          shot & $\|\mathbf{v}\|$ & $\|\mathbf{P}_3-\mathbf{q}\|$ \\
          \hline
          1 & 1         & $5.6\cdot10^{-1}$ \\
          2 & 1.5403024 & $3.0\cdot10^{-2}$ \\
          3 & 1.5707916 & $4.7\cdot10^{-6}$ \\
          4 & 1.5707963 & $3.4\cdot10^{-15}$ \\
          \hline
        \end{tabular}
        \end{center}
        In the fourth shot, $\|\mathbf{P}_3-\mathbf{q}\|<\delta$.
        We illustrate the points obtained on the fourth shooting in
        Figure~\ref{fig: sphere shooting}(b).
  \item The output is $(\mathbf{P}_0, \mathbf{P}_1, \mathbf{P}_2, \mathbf{P}_3)$ (see Figure~\ref{fig: sphere shooting}(b)) and $\|\mathbf{v}\| = 1.5707963\ldots$
\end{enumerate}

We note that the output agrees with $\pi/2$ up to $3\cdot10^{-14}$.

\begin{figure}[ht]
\centering
\begin{tikzpicture}[scale=2.2, >=stealth]
  \colorlet{cA}{orange!70!red}
  \colorlet{cB}{teal}
  \colorlet{vA}{cA!60!black}
  \colorlet{vB}{cB!60!black}
  \begin{scope}
    \shade[ball color=gray!10, opacity=0.5]
      (0,0) circle (1);
    \draw (0,0) circle (1);
    \draw[gray] (-1,0) arc (180:360:1 and 0.15);
    \draw[gray, dashed] (1,0)
      arc (0:180:1 and 0.15);
    \draw[gray, dashed] (0,0) -- (1,0);
    \draw[gray, dashed] (0,0) -- (0,1);
    \draw[gray, dashed] (0,0) -- (150:1);
    \draw[->] (1,0) -- (1.35,0)
      node[right] {$x_1$};
    \draw[->] (150:1) -- (150:1.35)
      node[above left=-2pt] {$x_2$};
    \draw[->] (0,1) -- (0,1.35)
      node[above] {$x_3$};
    \foreach \i/\j in
        {0/19.1, 19.1/38.2, 38.2/57.3} {
      \draw[very thick, cA]
        plot[domain=\i:\j, samples=40, smooth]
        (\x:{1-0.24*((\x-\i)/(\j-\i))^2});
      \draw[->, cA] (\j:0.8) -- (\j:0.955);
      \draw[thick, cA, fill=white]
        (\j:0.76) circle (0.035);
      \fill[cA] (\j:1) circle (0.04);
      \draw[->, thick, vA]
        (\i:1) -- ++({\i+90}:0.2);
    }
    \begin{scope}[cA]
      \node[left=1pt] at (19.1:0.76)
        {$\mathbf{x}_{\mathbf{v}}(t_1)$};
      \node[left=1pt] at (38.2:0.76)
        {$\mathbf{x}_{\mathbf{v}}(t_2)$};
      \node[above left=-1pt] at (57.3:0.76)
        {$\mathbf{x}_{\mathbf{v}}(t_3)$};
      \node[right=1pt] at (19.1:1)
        {$\mathbf{P}_1$};
      \node[above right=-1pt] at (38.2:1)
        {$\mathbf{P}_2$};
      \node[above right=-1pt] at (57.3:1)
        {$\mathbf{P}_3$};
      \node[right=1pt, vA] at (1,0.15)
        {$\mathbf{w}_0 = \mathbf{v}$};
      \node[right=1pt, vA] at (0.88,0.53)
        {$\mathbf{w}_1$};
      \node[vA] at (0.9,0.85)
        {$\mathbf{w}_2$};
    \end{scope}
    \draw[thick, fill=white]
      (90:1) circle (0.04);
    \node[above left=-1pt] at (90:1)
      {$\mathbf{q}$};
    \fill (1,0) circle (0.04);
    \node[below right] at (1,-0.06)
      {$\mathbf{P}_0 = \mathbf{p}$};
    \node at (0,-1.3)
      {(a) first shot, $\|\mathbf{v}\| = 1$};
  \end{scope}
  \begin{scope}[shift={(3.1,0)}]
    \shade[ball color=gray!10, opacity=0.5]
      (0,0) circle (1);
    \draw (0,0) circle (1);
    \draw[gray] (-1,0) arc (180:360:1 and 0.15);
    \draw[gray, dashed] (1,0)
      arc (0:180:1 and 0.15);
    \draw[gray, dashed] (0,0) -- (1,0);
    \draw[gray, dashed] (0,0) -- (0,1);
    \draw[gray, dashed] (0,0) -- (150:1);
    \draw[->] (1,0) -- (1.35,0)
      node[right] {$x_1$};
    \draw[->] (150:1) -- (150:1.35)
      node[above left=-2pt] {$x_2$};
    \draw[->] (0,1) -- (0,1.35)
      node[above] {$x_3$};
    \foreach \i/\j in {0/30, 30/60, 60/90} {
      \draw[very thick, cB]
        plot[domain=\i:\j, samples=40, smooth]
        (\x:{1-0.24*((\x-\i)/(\j-\i))^2});
      \draw[->, cB] (\j:0.8) -- (\j:0.955);
      \draw[thick, cB, fill=white]
        (\j:0.76) circle (0.035);
      \fill[cB] (\j:1) circle (0.04);
      \draw[->, thick, vB]
        (\i:1) -- ++({\i+90}:0.314);
    }
    \begin{scope}[cB]
      \node[below left=-1pt] at (30:0.76)
        {$\mathbf{x}_{\mathbf{v}}(t_1)$};
      \node[below left=-1pt] at (60:0.76)
        {$\mathbf{x}_{\mathbf{v}}(t_2)$};
      \node[left=1pt] at (90:0.76)
        {$\mathbf{x}_{\mathbf{v}}(t_3)$};
      \node[right=1pt] at (30:1)
        {$\mathbf{P}_1$};
      \node[above right=-1pt] at (60:1)
        {$\mathbf{P}_2$};
      \node[above left=-1pt] at (90:1)
        {$\mathbf{P}_3 = \mathbf{q}$};
      \node[right=1pt, vB] at (1,0.2)
        {$\mathbf{w}_0 = \mathbf{v}$};
      \node[right=1pt, vB] at (0.71,0.78)
        {$\mathbf{w}_1$};
      \node[above=1pt, vB] at (0.23,1.02)
        {$\mathbf{w}_2$};
    \end{scope}
    \fill (1,0) circle (0.04);
    \node[below right] at (1,-0.06)
      {$\mathbf{P}_0 = \mathbf{p}$};
    \node at (0,-1.3)
      {(b) last shot, $\|\mathbf{v}\| = \pi/2$};
  \end{scope}
\end{tikzpicture}
\caption{\normalsize The first and last shots of
Example~\ref{ex: sphere distance}.}
\label{fig: sphere shooting}
\end{figure}
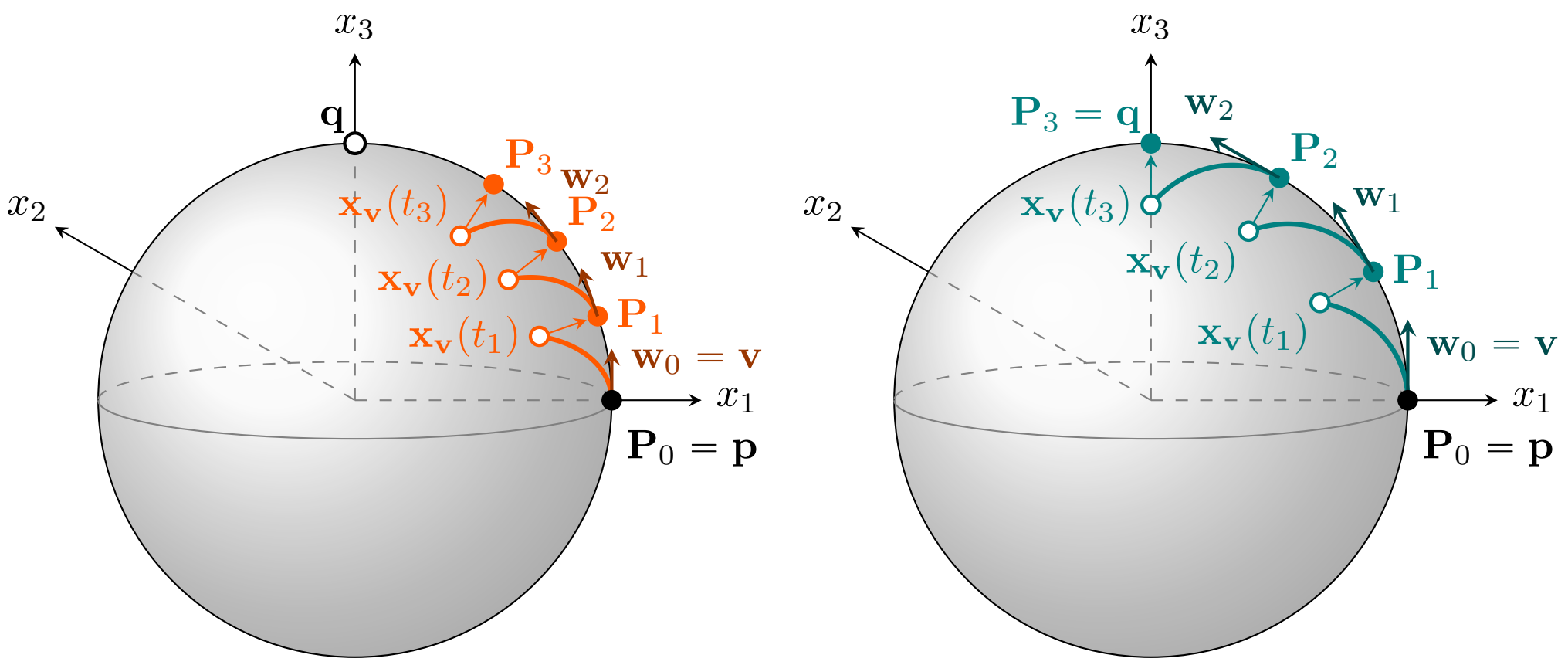

\end{example}

The output of Algorithm~\ref{alg:geodesic} is the length of a geodesic from $\mathbf{p}$ to $\mathbf{q}$, and hence only an upper bound
for $\mathrm{d}_X(\mathbf{p},\mathbf{q})$. To estimate the geodesic distance, it is
therefore useful to have a lower bound as
well. A lower bound comes from the ambient
space: every curve on $X$ is also a curve in
$\mathbb{R}^n$, so the geodesic distance on $X$
is at least the Euclidean distance between
$\mathbf{p}$ and $\mathbf{q}$, and for a
projective variety, it is at least the
distance in projective space.
In particular, if the
length computed by Algorithm~\ref{alg:geodesic} is close to the lower bound, then it is also close to $\mathrm{d}_X(\mathbf{p},\mathbf{q})$. We
state these bounds in the following lemmas.

\begin{lemma} \label{lem : lower bound}
    Let $X\subset \mathbb{K}^n$ be a smooth affine variety, $\mathbf{p},\mathbf{q}\in X$. Then
    \[ \dd_X(\mathbf{p},\mathbf{q})\geq \|\mathbf{p}-\mathbf{q}\|.\]
\end{lemma}
\begin{proof}
   This is due to the fact that the geodesics in $\mathbb{K}^n$ are lines.
\end{proof}

\begin{example}
For $\mathbf{p}$, $\mathbf{q}$ and $X$ from
Example~\ref{ex: sphere distance} we have
$\dd_X(\mathbf{p},\mathbf{q})\approx
1.571$, while 
$\|\mathbf{p}-\mathbf{q}\| = \sqrt{2}\approx
1.414$.
\end{example}

\begin{lemma}
    Let $X$ be a projective variety, $[\mathbf{p}],[\mathbf{q}]\in X$. Then
    \begin{align*}
        \mathrm{d}_X([\mathbf{p}],[\mathbf{q}])\geq\arccos{\frac{|\langle \mathbf{p},\mathbf{q}\rangle|}{\|\mathbf{p}\|\|\mathbf{q}\|}}.
    \end{align*}
\end{lemma}
\begin{proof}
    The right-hand side of the inequality is the expression for the geodesic distance on $\mathbb{KP}^n$.
\end{proof}

\subsection{Implementation and benchmarks}

We implemented Algorithm~\ref{alg:geodesic} in
\texttt{Julia}, using \textsc{Oscar}~\cite{OSCAR} to
compute the geodesic equations in
line~\ref{line:equations} symbolically, as
rational functions in $\mathbf{x}$ and
$\mathbf{x}'$, and
OrdinaryDiffEq~\cite{DifferentialEquationsJl}
for the Runge--Kutta method.

The code is
available at
\begin{center}
{\color{teal}\url{https://github.com/ymukhina/geodesics-on-algebraic-varieties}.}
\end{center}

We compare Algorithm~\ref{alg:geodesic} with
the parametric case of
Remark~\ref{rem: variants} on the following
varieties, with $K = 10$,
$\varepsilon = 10^{-6}$,
$\delta = 10^{-11}$ and $N=20$. For each of them, below we give
the defining polynomials $\mathbf{F}$ used in
Algorithm~\ref{alg:geodesic} and the
parametrization $\varphi$ used in the
parametric case. 
We note that the Cayley cubic and the
cone over $\operatorname{Gr}(2,4)$ are singular
at the origin. Nevertheless, Algorithm~\ref{alg:geodesic}
works as long as the computed curve stays in
the smooth part of $X$ and points
$\mathbf{p}$ and $\mathbf{q}$ are smooth.

\begin{enumerate}
  \item \emph{The sphere} $X\subset\R^3$:
    \[
      \mathbf{F} = (x_1^2 + x_2^2 + x_3^2 - 1), \qquad
      \varphi(s) = \frac{1}{1 + s_1^2 + s_2^2}\,
      \bigl(2s_1,\ 2s_2,
      \ s_1^2 + s_2^2 - 1\bigr).
    \]
  \item \emph{The triaxial ellipsoid}
    $X\subset\R^3$:
    \[
      \mathbf{F} = (4x_1^2 + 9x_2^2 + 36x_3^2 - 36),
      \qquad
      \varphi(s) = \frac{1}{1 + s_1^2 + s_2^2}\,
      \bigl(6s_1,\ 4s_2,
      \ s_1^2 + s_2^2 - 1\bigr).
    \]
  \item \emph{The hyperbolic paraboloid}
    $X\subset\R^3$:
    \[
      \mathbf{F} = (x_1^2 - x_2^2 - 2x_3), \qquad
      \varphi(s) = \bigl(s_1,\ s_2,\
      \tfrac{1}{2}(s_1^2 - s_2^2)\bigr).
    \]
  \item \emph{The Cayley cubic} $X\subset\R^3$:
    \[
      \mathbf{F} = (x_1x_2x_3 + x_1x_2 + x_1x_3 + x_2x_3),
      \qquad
      \varphi(s) = \Bigl(\frac{S}{s_1},\
      \frac{S}{s_2},\ S\Bigr),
      \quad S = -(s_1 + s_2 + 1).
    \]
  \item \emph{The rotation group}
    $\mathrm{SO}(3)\subset\R^{3\times3} = \R^9$,
    parametrized by the Cayley transform. We
    write a point of $\R^9$ as the matrix
    \[
      A = \begin{pmatrix}
        x_1 & x_2 & x_3\\
        x_4 & x_5 & x_6\\
        x_7 & x_8 & x_9
      \end{pmatrix},
    \]
    and $\mathbf{F} = (F_1, \ldots, F_6)$ consists of the six entries of
    $A^TA - I$ on and above the diagonal, that
    is, the columns of $A$ are orthonormal:
    \[
      \begin{aligned}
        &x_1^2 + x_4^2 + x_7^2 - 1, &
        &x_1x_2 + x_4x_5 + x_7x_8,\\
        &x_2^2 + x_5^2 + x_8^2 - 1, &
        &x_1x_3 + x_4x_6 + x_7x_9,\\
        &x_3^2 + x_6^2 + x_9^2 - 1, &
        &x_2x_3 + x_5x_6 + x_8x_9.
      \end{aligned}
    \]
    The parametrization is
    \[
      \varphi(s) = \frac{1}{1 + s_1^2 + s_2^2
        + s_3^2}
      \begin{pmatrix}
        1 - s_1^2 - s_2^2 + s_3^2
          & -2(s_1 + s_2s_3)
          & 2(s_1s_3 - s_2)\\
        2(s_1 - s_2s_3)
          & 1 - s_1^2 + s_2^2 - s_3^2
          & -2(s_1s_2 + s_3)\\
        2(s_2 + s_1s_3)
          & 2(s_3 - s_1s_2)
          & 1 + s_1^2 - s_2^2 - s_3^2
      \end{pmatrix}.
    \]
    Then $\mathcal{V}(\mathbf{F}) = \mathrm{O}(3)$, and
    $\mathrm{SO}(3)$ is the connected component
    of $\mathcal{V}(\mathbf{F})$ that contains the identity
    matrix. We leave out the equation
    $\det A = 1$: a curve never leaves the
    connected component of its initial point,
    and Algorithm~\ref{alg:geodesic} needs
    exactly $\operatorname{codim} X = 6$
    equations.
  \item \emph{The cone over the Grassmannian}
    $\operatorname{Gr}(2,4)$, $X\subset\R^6$:
    \[
      \mathbf{F} = (x_1x_6 - x_2x_5 + x_3x_4), \qquad
      \varphi(s) = \bigl(s_1^2,\ s_1s_4,\
      s_1s_5,\ -s_1s_2,\ -s_1s_3,\
      s_2s_5 - s_3s_4\bigr).
    \]
\end{enumerate}

Table~\ref{tab:benchmarks} shows the results
for several pairs of points
$\mathbf{p} = \varphi(\mathbf{s}_\mathbf{p})$ and
$\mathbf{q} = \varphi(\mathbf{s}_\mathbf{q})$: the
lower bound of Lemma~\ref{lem : lower bound}, the lengths of geodesics and
running times. The times include the computation
of the geodesic equations. All computations
were done on an Apple M2 Pro with 32\,GB of
memory, with \texttt{Julia}~1.12.5,
\textsc{Oscar}~1.8.2 and OrdinaryDiffEq~7.8.1.

\begin{table}[ht]
\centering
\setlength{\tabcolsep}{5pt}
\begin{tabular}{@{}lllrrrrr@{}}
\toprule
& & & & \multicolumn{2}{c}{implicit}
  & \multicolumn{2}{c@{}}{parametric} \\
\cmidrule(lr){5-6} \cmidrule(l){7-8}
variety & $\mathbf{s}_\mathbf{p}$ & $\mathbf{s}_\mathbf{q}$ & lower
  & length & time & length & time \\
\midrule
\multirow{3}{2.1cm}{\raggedright%
\textbf{sphere}}
  & $(0, 0.5)$
  & $(0, 2)$
  & 1.20
  & 1.29 & $0$
  & 1.29 & $0$ \\
  & $(0.3, 0.4)$
  & $(1, 2)$
  & 1.28
  & 1.38 & $0$
  & 1.38 & 8.6 \\
  & $(0.6, 0.8)$
  & $(-0.8, 0.6)$
  & 1.41
  & 1.57 & $0$
  & 1.57 & $0$ \\
\midrule
\multirow{3}{2.1cm}{\raggedright%
\textbf{triaxial ellipsoid}}
  & $(0, 0.5)$
  & $(0, 2)$
  & 1.20
  & 1.50 & $0$
  & 1.50 & $0$ \\
  & $(0.3, 0.4)$
  & $(1, 2)$
  & 1.34
  & 1.74 & $0$
  & 1.74 & 1.6 \\
  & $(0.5, -0.5)$
  & $(-0.3, 0.2)$
  & 4.16
  & 4.29 & $0$
  & 4.29 & $0$ \\
\midrule
\multirow{3}{2.1cm}{\raggedright%
\textbf{hyperbolic paraboloid}}
  & $(0, -1)$
  & $(0, 2)$
  & 3.35
  & 4.11 & $0$
  & 4.11 & $0$ \\
  & $(1, 0)$
  & $(2, 0)$
  & 1.80
  & 1.81 & $0$
  & 1.81 & $0$ \\
  & $(1, 2)$
  & $(-1, 0.5)$
  & 3.13
  & 3.14 & $0$
  & 3.14 & $0$ \\
\midrule
\multirow{3}{2.1cm}{\raggedright%
\textbf{Cayley cubic}}
  & $(2, 3)$
  & $(3, 4)$
  & 2.03
  & 2.03 & $0$
  & 2.03 & $0$ \\
  & $(-3, -2)$
  & $(-7, -3)$
  & 5.10
  & 5.10 & $0$
  & 5.10 & $0$ \\
  & $(1, 0.1)$
  & $(0.1, 1)$
  & 26.73
  & 26.95 & $0$
  & 26.95 & $0$ \\
\midrule
\multirow{6}{2.1cm}{\raggedright%
\textbf{\boldmath$\mathrm{SO}(3)$}}
  & $(1, 2, 3)$
  & $(4, 5, 6)$
  & 0.74
  & 0.75 & 17
  & 0.75 & $0$ \\
  & $(-5, 3, 4)$
  & $(-1, 7, 100)$
  & 2.27
  & 2.64 & 35
  & 2.64 & 0.2 \\
  & $(0, 1, 0)$
  & $(3, 0, 5)$
  & 2.81
  & 4.10 & 49
  & 4.10 & 0.1 \\
  & $(1, -2, 3)$
  & $(4, 5, -6)$
  & 2.09
  & {\color{teal}2.36} & 24
  & {\color{orange!70!red}6.53} & 0.1 \\
  & $(10, -2, 3)$
  & $(-32, 17, -15)$
  & 0.88
  & {\color{teal}0.90} & 17
  & {\color{orange!70!red}7.99} & 0.3 \\
  & $(-5, 3, 4)$
  & $(-1, 7, -24)$
  & 2.61
  & {\color{teal}3.32} & 35
  & {\color{orange!70!red}5.56} & 0.1 \\
\midrule
\multirow{3}{2.1cm}{\raggedright%
\textbf{cone over \boldmath$\mathrm{Gr}(2,4)$}}
  & $(1, 2, 3, 4, 5)$
  & $(6, 7, 8, 9, 10)$
  & 101.86
  & 101.86 & $0$
  & 101.86 & $0$ \\
  & $(1, -2, 3, 4, 5)$
  & $(6, 7, 8, 9, 10)$
  & 105.41
  & 105.68 & $0$
  & 105.68 & $0$ \\
  & $(0.1, 2, 1, 2.2, 3)$
  & $(1, 0.2, 1.3, 2, 2)$
  & 6.67
  & 6.76 & $0$
  & 6.76 & $0$ \\
\bottomrule
\end{tabular}
\caption{\normalsize Lengths and running times
of the geodesics computed by
Algorithm~\ref{alg:geodesic} in the implicit
and in the parametric cases with the lower
bound of Lemma~\ref{lem : lower bound}.
Lengths are rounded to two decimals. Runtimes
are in seconds (runtimes below $0.1$ seconds
are shown as $0$).}
\label{tab:benchmarks}
\end{table}

Depending on the example, either the implicit
or the parametric version is faster. Both
versions find the same geodesic in all
examples except three pairs on
$\mathrm{SO}(3)$, colored in
Table~\ref{tab:benchmarks}.
There, the parametric version finds a longer
geodesic, because $\varphi$ does not cover all
of $\mathrm{SO}(3)$: it misses the rotations by
the angle $\pi$. 
In these three pairs, the
shortest geodesic from $\mathbf{p}$ to
$\mathbf{q}$ passes through such a rotation,
so the parametric version cannot follow it.

We note that on the unit sphere and on
$\mathrm{SO}(3)$, the geodesic distance is
given by a closed formula. For the sphere, the
formula is
$\arccos\langle\mathbf{p},\mathbf{q}\rangle$~%
\cite{Lee2018}, and on $\mathrm{SO}(3)$ it is
$\sqrt{2}\,\theta$, where $\theta$ is a certain
rotation angle (see \cite{Moakher2002}).
Using these formulas, one may check the
results of our implementation as well. For the
chosen pairs of points, we note that the
lengths computed by
Algorithm~\ref{alg:geodesic} agree with the
exact distance up to a relative error of less
than $10^{-12}$.

We do not compare our implementation with other
software: packages such as
Manopt~\cite{Manopt},
Geomstats~\cite{Geomstats} and
Manifolds.jl~\cite{ManifoldsJl} treat specific
manifolds, such as spheres and rotation groups,
often by closed formulas like the ones above,
and to our knowledge they do not provide
geodesics on our other four considered varieties.

\begin{figure}[ht]
  \centering
\includegraphics[width=0.45\textwidth]
{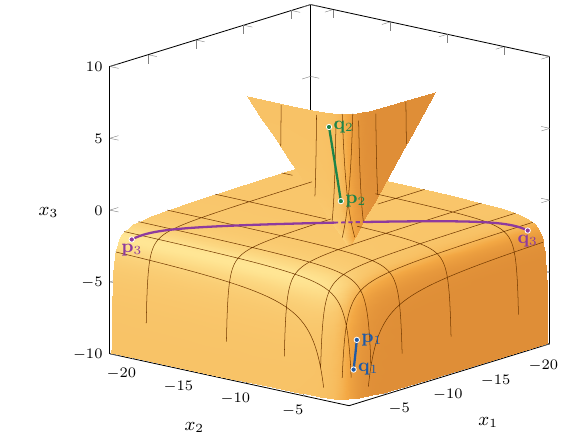}
  \caption{Geodesics on the Cayley cubic
  between chosen points $\mathbf{p}_i$ and
  $\mathbf{q}_i$.}
  \label{fig:cayley geodesics}
\end{figure}

We finish this section with a picture of
geodesics on the Cayley cubic. For the points
\begin{alignat*}{2}
  \mathbf{p}_1&=\varphi(2,3),
  &\qquad \mathbf{q}_1&=\varphi(3,4),\\
  \mathbf{p}_2&=\varphi(-3,-2),
  &\qquad \mathbf{q}_2&=\varphi(-7,-3),\\
  \mathbf{p}_3&=\varphi(1,0.1),
  &\qquad \mathbf{q}_3&=\varphi(0.1,1),
\end{alignat*}
from Table~\ref{tab:benchmarks}, we obtain the
geodesics from $\mathbf{p}_i$ to
$\mathbf{q}_i$ shown in
Figure~\ref{fig:cayley geodesics}. To draw
them, we ran Algorithm~\ref{alg:geodesic} with
$K=100$, $\varepsilon=10^{-6}$,
$\delta=10^{-11}$ and $N=20$.

\subsection*{Acknowledgements}
We thank Bernd Sturmfels and Gleb Pogudin for their suggestions and guidance, and thank Emil Horobet for his support of and enthusiasm for our project. We express our gratitude to the Max Planck Institute for Mathematics in the Sciences for providing accommodation for PNP and YM during the end stages of the work. CdK was supported by the European Research Council through the synergy grant UNIVERSE+, 101118787. 
YM was supported by the French ANR-22-CE48-0008 OCCAM project and would like
to gratefully acknowledge this support.
This work has been partially funded by the European Union under the Grant Agreement no. 101044561, POSALG. 
PNP acknowledges partial support by the Studienstiftung des deutschen Volkes and the Methusalem grant METH/21/03 - long term structural funding of the Flemish Government.

\noindent {\tiny Views and opinions
expressed are however those of the authors only and do not necessarily reflect those of the European Union or the European
Research Council Executive Agency. Neither the European Union nor the granting authority can be held responsible for them.}

\appendix

\section{Implicit geodesic equations on projective varieties}
\label{appendix projective}
In this section we will discuss how to arrive at the implicit geodesic equations for projective varieties stated in Theorem~\ref{impl form proj : thm}. To this end, we will use the geodesic equation in the form $\nabla^X_{\gamma'}\gamma'=0$, where $\nabla^X$ denotes the Levi-Civita connection on the given projective variety $X$ with metric induced by  the Fubini--Study metric.

We will need two classical results from the literature. The first concerns embedded submanifolds $N$ of a Riemannian manifold $(M,g)$ equipped with the induced metric. 
\begin{theorem}[Gauss Formula, see e.g.\ Theorem 8.2 in~\cite{Lee2018}]
    \label{second fundamental form : thm}
    Let $(M,g)$ be a Riemannian manifold, $N$ an embedded submanifold with the induced metric and $V,W\in\X(N)$ vector fields on $N$. Then
    \begin{align*}
        \nabla^M_V W=\nabla^N_V W +\mathrm{I\!I}(V,W), 
    \end{align*}
    where $\mathrm{I\!I}(V,W)$ is the \emph{second fundamental form}, which takes values in $TN^{\perp}\subset TM|_N$. 
\end{theorem}

The second statement that we will need is a similar formula for \emph{Riemannian submersions}. These are submersions $\pi\colon(\tilde{M},\tilde{g})\rightarrow (M,g)$ satisfying
$$\tilde{g}(V,W)=g(d\pi (V),d\pi (W))\quad \forall\, V,W\in (\ker d\pi)^\perp\subset T\tilde{M}.$$
In particular, the projection from the sphere with the round metric to projective space with the Fubini--Study metric is a Riemannian submersion.

Given a Riemannian submersion $\pi\colon(\tilde{M},\tilde{g})\rightarrow (M,g)$, the tangent bundle of $\tilde{M}$ decomposes into the kernel of $d\pi$, called \emph{vertical space} $\mathcal{V}$ and its orthogonal complement with respect to $\tilde{g}$, called \emph{horizontal space} $\mathcal{H}$. Every vector field $V\in \X(M)$ has a unique \emph{horizontal lift} $\tilde{V}\in \X(\tilde{M})$, s.t.\ $\tilde{V}_x\in\mathcal{H}_x$ for all $x\in\tilde{M}$ and~$d\pi(\tilde{V})=V$.

\begin{proposition}
[See e.g.\ Proposition 4.5.1 in~\cite{Petersen2016}]
\label{riemannian submersions : prop}
    Let $(\tilde{M},\tilde{g})\to (M,g)$ be a Riemannian submersion. We denote the Levi-Civita connections of $\tilde{g}$ and $g$ by $\tilde{\nabla}$ and $\nabla$, respectively. Then for all $V,W\in\X(M)$, we have
    \begin{align*}
        \tilde{\nabla}_{\tilde{V}}\tilde{W}=\widetilde{\nabla_V W}+\tfrac{1}{2}[\tilde{V},\tilde{W}]^{\mathcal{V}},
    \end{align*}
    where $\cdot^{\mathcal{V}}$ denotes the vertical part of a vector field on $\tilde{M}$.
\end{proposition}

In the following, we will sometimes identify $\mathbb{C}^{n+1}\cong \mathbb{R}^{2n+2}$ via $(x_0+\ii y_0,\dots,x_n+\ii y_n)\mapsto (x_0,y_0,\dots,x_n,y_n)$. Under this identification, the multiplication by $\ii$ is given by the matrix
$$J=\begin{psmallmatrix}
    0 & -1   & & &\\
    1 & 0 & & &\\
     & & \ddots & & \\
     & &  & 0 & -1 \\
     & & & 1 & 0 
\end{psmallmatrix}.$$
Moreover, the vanishing of the polynomials is equivalent to the vanishing of their real and imaginary parts. When working with the identification $\C^{n+1}\cong\R^{2n+2}$ we will therefore consider the Jacobian 
\begin{align}
   \mathcal{J}_{\mathbf{F}}=\begin{pmatrix}
       \del\operatorname{Re}(F_1)/\del x_0 &  \del\operatorname{Re}(F_1)/\del y_0 &\cdots & \del\operatorname{Re}(F_1)/\del x_n & \del\operatorname{Re}(F_1)/\del y_n\\
       \del\operatorname{Im}(F_1)/\del x_0 &\del\operatorname{Im}(F_1)/\del y_0  & \cdots & \del\operatorname{Im}(F_1)/\del x_n&\del\operatorname{Im}(F_1)/\del y_n\\
       \vdots & \vdots&\ddots & \vdots &\vdots\\
       \del\operatorname{Re}(F_k)/\del x_0 &  \del\operatorname{Re}(F_k)/\del y_0 &\cdots & \del\operatorname{Re}(F_k)/\del x_n & \del\operatorname{Re}(F_k)/\del y_n\\
       \del\operatorname{Im}(F_k)/\del x_0 &\del\operatorname{Im}(F_k)/\del y_0  & \cdots & \del\operatorname{Im}(F_k)/\del x_n&\del\operatorname{Im}(F_k)/\del y_n
   \end{pmatrix}.  
\end{align}
Recall furthermore the definition of the Euler vector field $E(\mathbf{x})=\sum_i x_i\del_{x_i}$ on $\R^{n+1}$ (or, for $\mathbb{K}=\C$, $E(\mathbf{x},\mathbf{y})=\sum_i x_i\del_{x_i}+y_i\del_{y_i}$ on $\R^{2n+2}$) and of the complex Euler vector field $Z(\mathbf{z})=\sum_i z_i\del_{z_i}$ on $\C^{n+1}\cong \R^{2n+2}$, where $z_i=x_i+\ii y_i$ and $\del_{z_i}=\tfrac{1}{2}(\del_{x_i}-\ii\del_{y_i})$. Using Theorem~\ref{second fundamental form : thm} and Proposition~\ref{riemannian submersions : prop}, we prove the following:
\begin{proposition}\label{appendix projective impl : prop}
    Let $X=\mathcal{V}(F_1,\dots,F_k)\subset \mathbb{KP}^n$ be a smooth projective variety and   let $\gamma\colon [a,b]\rightarrow X$ be a curve on $X$. We consider its horizontal lift to\footnote{Here $N=n$ if $\mathbb{K}=\R$ and $N=2n+1$ if $\mathbb{K}=\C$.} $S^N\subset \mathbb{K}^{n+1}$ as a curve in $\mathbb{K}^{n+1}$. \begin{enumerate}
        \item If $X$ is real projective, then $\gamma$ is a geodesic if and only if
        $$\tilde{\gamma}''\in\operatorname{span}_{\R}\{E, \operatorname{grad} F_1,\dots,\operatorname{grad} F_k\},$$
        where $\operatorname{grad}F_i=\sum_k\tfrac{\del F_i}{\del x_k}\del_{x_k}$ and for $\tilde{\gamma}(t)=(x_0(t),\dots,x_n(t))$ we write $\tilde{\gamma}''=\sum_i x_i''\del_{x_i}$.
        \item If $X$ is complex projective, then $\gamma$ is a geodesic if and only if
        $$\tilde{\gamma}''\in\operatorname{span}_{\C}\{Z, \overline{\operatorname{grad} F_1},\dots,\overline{\operatorname{grad} F_k}\},$$
        where $\overline{\operatorname{grad} F_i}=\sum_k \overline{\frac{\del F_i}{\del z_k}}\del_{z_k}$ and for $\tilde{\gamma}(t)=(z_0(t),\dots,z_n(t))$ we write $\tilde{\gamma}''=\sum_i z_i''\del_{z_i}$. 
    \end{enumerate}
\end{proposition}
\begin{proof}
Let $X=\mathcal{V}(F_1,\dots,F_k)\subset \mathbb{KP}^n$. First, note that by Theorem~\ref{second fundamental form : thm} a curve $\gamma\subset X$ is a geodesic, i.e.\ $\nabla^X_{\gamma'}\gamma'=0$, if and only if $\nabla^{\mathbb{KP}^n}_{\gamma'}\gamma'\in TX^{\perp}$. Moreover, observe that $TX=d\tilde{\pi}(\ker \mathcal{J}_\mathbf{F})$, where $\tilde{\pi}\colon \mathbb{K}^{n+1}\to \mathbb{KP}^n$ and $\mathbf{F}=(F_1,\dots,F_k)$ if $\mathbb{K}=\R$ and $\mathbf{F}=(\operatorname{Re}(F_1),\operatorname{Im}(F_1),\dots,\operatorname{Re}(F_k),\operatorname{Im}(F_k))$ if $\mathbb{K}=\C$, respectively. Since $\tilde{\pi}=\pi\circ\mathrm{pr}_{S^N}$, where $\operatorname{pr}_{S^N}\colon \mathbb{K}^{n+1}\setminus\{0\}\to S^N$, $\mathbf{x}\mapsto \tfrac{\mathbf{x}}{\|\mathbf{x}\|}$, we find that $\gamma$ is a geodesic if and only if 
$$\nabla^{\mathbb{KP}^n}_{\gamma'}\gamma'\in d\pi \circ d\mathrm{pr}_{S^N}(\ker \mathcal{J}_\mathbf{F})^{\perp}.$$

By Proposition~\ref{riemannian submersions : prop}, we have
\begin{align*}
    \nabla^{S^N}_{\tilde{\gamma}'}\tilde{\gamma}'=\widetilde{\nabla^{\mathbb{KP}^n}_{\gamma'}\gamma'}+\tfrac{1}{2}[\tilde{\gamma}',\tilde{\gamma}']^{\mathcal{V}}=\widetilde{\nabla^{\mathbb{KP}^n}_{\gamma'}\gamma'}.
\end{align*}
Therefore, $\gamma$ is a geodesic if and only if
$$\nabla^{S^N}_{\tilde{\gamma}'}\tilde{\gamma}'\in d\operatorname{pr}_{S^N}(\ker \mathcal{J}_\mathbf{F})^{\perp}.$$
Moreover, it follows from Theorem~\ref{second fundamental form : thm} and the fact that for the Euclidean connection $\nabla^{\mathrm{E}}$, we have $\nabla^{\mathrm{E}}_{\tilde{\gamma}'}\tilde{\gamma}'=\tilde{\gamma}''$ that
$$\nabla^{S^N}_{\tilde{\gamma}'}\tilde{\gamma}'=\tilde{\gamma}''+\mathrm{I\!I}(\tilde{\gamma}',\tilde{\gamma}').$$
Here, $\mathrm{I\!I}(\tilde{\gamma}',\tilde{\gamma}')\in (TS^N)^\perp=\operatorname{span}\{E\}$.

Putting everything together (and using that $\ker d\operatorname{pr}_{S^N}=\operatorname{span}_{\R}\{E\}$), we find that $\gamma$ is a geodesic if and only if
\begin{align}\label{eq: appendix in proof of proposition}
    \tilde{\gamma}''\in \operatorname{span}_{\R}\{E\}\oplus (\ker \mathcal{J}_{\mathbf{F}})^{\perp}.
\end{align}
If $\mathbb{K}=\mathbb{R}$, this finishes the proof. If $\mathbb{K}=\mathbb{C}$, observe that $$Z(\mathbf{z})=\tfrac{1}{2}\left(E(\mathbf{z})-\ii JE(\mathbf{z})\right).$$
Using this and the identification $\C^n\cong\R^{2n}$, we find that \eqref{eq: appendix in proof of proposition} is equivalent to
$$\tilde{\gamma}''\in\operatorname{span}_{\C}\{Z\}\oplus(\ker \mathcal{J}_{\mathbf{F}})^{\perp},\quad 0=\operatorname{Im}\langle\tilde{\gamma}'',Z\rangle=\operatorname{Im}\langle\tilde{\gamma}'',\tilde{\gamma}\rangle,$$
where now $\mathcal{J}_{\mathbf{F}}$ denotes the complex Jacobian $\mathcal{J}_{\mathbf{F},ij}=\del F_i/\del z_j$. Note that the condition $\operatorname{Im}\langle\tilde{\gamma}'',\tilde{\gamma}\rangle=0$ is always satisfied for horizontal curves since
$$\langle \tilde{\gamma}',\tilde{\gamma}\rangle=0\implies \langle\tilde{\gamma}'',\tilde{\gamma}\rangle+\langle \tilde{\gamma}',\tilde{\gamma}'\rangle=0.$$ Moreover, 
$$(\ker \mathcal{J}_{\mathbf{F}})^{\perp}=\operatorname{span}\{\overline{\operatorname{grad} F_i}\},$$
where $\overline{\operatorname{grad} F_i}=\sum_k \overline{\tfrac{\del F_i}{\del z_k}}\del_{z_k}$. 
\end{proof}
The following lemma shows that appending $E$ or $Z$ as an extra row to $\overline{\mathcal{J}_{\mathbf{F}}}$ indeed increases the rank.
\begin{lemma}[See also Lemma 2.8 in \cite{Draisma2015}]\label{lem: appendix for rank of augm jac}
If $X=\mathcal{V}(\mathbf{F})$ is a real projective variety, then $$E(\mathbf{x})\in (\operatorname{span}_{\R}\{\operatorname{grad}F_i(\mathbf{x})\})^{\perp}\quad\forall \mathbf{x}\in X.$$  
Similarly, if $X=\mathcal{V}(\mathbf{F})$ is a complex projective variety, then $$Z(\mathbf{x})\in (\operatorname{span}_{\C}\{\overline{\operatorname{grad}F_i(\mathbf{x})}\})^{\perp}\quad\forall \mathbf{x}\in X.$$
\end{lemma}
\begin{proof}
    We will only provide a short proof for the complex case. The real case works similarly. Recall that for a homogeneous polynomial $F$ of degree $r$ the (complex) Euler vector field satisfies $Z[F]=dF(Z)=r\cdot F$. In particular, for every $\mathbf{x}\in \mathcal{V}(F)$, we have $dF(\mathbf{x})Z(\mathbf{x})=0$ and hence $\langle Z(\mathbf{x}),\overline{\operatorname{grad}F(\mathbf{x})}\rangle=0$.
\end{proof}
\small
\bibliographystyle{abbrv}
\bibliography{ref}
\end{document}